\documentclass[12pt]{article}
 \usepackage{amsmath}
 \usepackage{amssymb}
 \usepackage[noadjust]{cite}
 \usepackage{color}
 \usepackage{curves}
 \usepackage{enumerate}
 \usepackage{epic}
 \usepackage{epsfig}
 \usepackage{float}
 \usepackage{graphics}
 \usepackage{graphicx}
 \usepackage{latexsym}
 \usepackage{multicol}
 \usepackage{subfigure}
 \usepackage{theorem}

 \newtheorem{thm}{Theorem}[section]
 \newtheorem{hyp}{Hypothesis}
 \newtheorem{cor}[thm]{Corollary}

 \theorembodyfont{\upshape}
 \newtheorem{rem}[thm]{Remark}
 \numberwithin{equation}{section}

 \DeclareMathOperator{\ccc}{c.c.}

 \DeclareMathOperator{\DD}{D}

 \DeclareMathOperator{\imag}{Im}
 
 \DeclareMathOperator{\Lip}{Lip}

 \DeclareMathOperator{\real}{Re}

 \DeclareMathOperator{\GL}{GL}

 \newcommand{\GC}{\mathbb{C}}

 \newcommand{\GN}{\mathbb{N}}

 \newcommand{\GR}{\mathbb{R}}
 \newcommand{\GS}{\mathbb{S}}
 \newcommand{\GSE}{\mathbb{SE}}
 \newcommand{\GSO}{\mathbb{SO}}
 
 \newcommand{\GZ}{\mathbb{Z}}
  \newcommand{\proof}{\textbf{Proof:\quad}}

 \newcommand{\ov}{\overline}

 \newcommand{\sdp}{\dot{+}}

 \newcommand{\newl}{\newline\newline}

  \DeclareRobustCommand{\qed}{%
  \ifmmode \mathqed
  \else
    \leavevmode\unskip\penalty9999 \hbox{}\nobreak\hfill
    \quad\hbox{$\Box$\normalsize}%
  \fi
}
\begin{document}

\title{Epicyclic drifting in anisotropic excitable media with multiple
inhomogeneities}
\author{P. Boily, \small \texttt{pboily@uottawa.ca} \\ \small Department of Mathematics and Statistics, University of
Ottawa\\
\small Ottawa K1N 6N5, Canada} \maketitle \abstract{Spirals have
been studied from a dynamical system perspective starting with
Barkley's seminal papers linking a wide class of spiral wave
dynamics to the Euclidean symmetry of the excitable media in which
they are observed. However, in order to explain certain
non-Euclidean phenomena, such as anchoring and epicyclic drifting,
LeBlanc and Wulff introduced a single translational
symmetry-breaking perturbation to the center bundle equation and
showed that rotating waves may be attracted to a non-trivial
solution manifold and travel epicyclically around the perturbation
center.
\par In this paper, we continue the (model-independent) investigation of
the effects of inhomogeneities on spiral wave dynamics by studying
epicyclic drifting in the presence of: a) $n$ simultaneous
translational symmetry-breaking terms, with $n > 1$, and b) a
combination of a single rotational symmetry-breaking term and a
single translational symmetry-breaking term. These types of forced
Euclidean symmetry-breaking provide a much more realistic model of
certain excitable media such as cardiac tissue. However, the main
theoretical tool used by LeBlanc and Wulff can only be applied to
their particular perturbation: we show how an averaging theorem of
Hale can be modified to analyze our two more general scenarios and
state the conditions under which epicyclic drifting takes place in
the general case. In the process, we recover LeBlanc and Wulff's
specific result. Finally, we illustrate our results with the help of
a simple numerical simulation of a modified bidomain
model.}\newpage\noindent \textbf{Keywords}: symmetry-breaking,
integral manifold, epicycle, spiral wave, excitable media,
averaging, center bundle equation, center manifold reduction
theorem.

\normalsize
\section{Introduction}
Spiral are found in numerous excitable media
\cite{LOPS,GZM,ZM,MPMPV,WR,YP,Detal,R1,Wetal,J,MAK,B1,B2,BKT} and
they give rise to beautiful imagery. While this in itself might
yield enough interest to study them, there is also (at least) one
serious reason to do so: spiral waves have been linked to cardiac
arrhythmias (to disruptions of the heart's normal electrical cycle)
\cite{W,Wetal,Detal,KS}. Furthermore, \begin{quote}most arrhythmias
are harmless but if they are ``re-entrant in nature and [...] occur
 because of the spatial distribution of cardiac tissue'' they can
seriously hamper the pumping mechanism of the heart and lead to
death \cite[p. 401]{KS}.\end{quote} As a result, a fuller
understanding of spiral wave dynamics in these media becomes
imperative. \subsubsection*{The equivariant dynamical system
approach} In recent years, one of the most rewarding approach to the
study of spiral waves is based on Barkley's initial observation that
the observed transition from rotating to modulated rotating wave can
be explained \textit{via} a Hopf bifurcation together with the
underlying Euclidean symmetries of the governing reaction-diffusion
equations \cite{B1,B2} (\textit{i.e.}: the semi-flow generated by
the dynamical system commutes with the
\begin{equation}
u(t,x)\longmapsto
u(t,x_1\cos\,\theta-x_2\sin\,\theta+p_1,x_1\sin\,\theta+
x_2\cos\,\theta+p_2), \label{uaction}
\end{equation}
where $(\theta,p_1,p_2)\in\mathbb{S}^1\times\mathbb{R}^2\simeq
\GSE(2)$ and $x\in \GR^2$ \cite{Wulff,DMcK}). This lead Barkley to
formulate a simpler \textit{ad hoc} 5-dimensional ODE system with
Euclidean symmetry replicating the above transition \cite{BK}. \par
Sandstede, Scheel and Wulff then proved a general center manifold
reduction theorem (CMRT) for relative equilibria and relative
periodic solutions in spatially extended infinite-dimensional
Euclidean-equivariant dynamical systems, providing a mathematical
justification of Barkley's insight \cite{SSW1,SSW2,SSW3,SSW4,FSSW}.
{However, this center manifold reduction theorem requires that the
spiral wave satisfy certain spectral gap conditions, which often
fail \cite{S}. Sandstede and Scheel have developed a comprehensive
theory of spiral instabilities using techniques of spatial dynamics
\cite{abscon2,abscon1} to deal with such a situation.}\par Other
methods are also used to reduce the dynamics to
finite-di\-men\-sio\-nal systems (such as the kinematic model using
the curvature of the wave as a driving mechanism \cite{MDZ}), but
the equivariant dynamical system approach has the advantage that it
can often provide universal, model-independent explanations and
predictions regarding the dynamics and bifurcations of spiral waves.
For example, the fore-mentioned `Hopf bifurcation' from rigid
rotation to quasi-periodic meandering has been observed in
numerically \cite{BKT} and experimentally \cite{LOPS}. Another
example is provided by the anchoring/repelling of spiral waves
on/from a site of inhomogeneity, which has been observed in
numerical integrations of an Oregonator system \cite{MPMPV}, in
photo-sensitive chemical reactions \cite{ZM} and in cardiac tissue
\cite{Detal}: using a model-independent approach based on forced
symmetry-breaking, LeBlanc and Wulff showed that anchoring/repelling
of rotating waves is a generic property of systems in which the
translation symmetry is broken by a small perturbation
\cite{LW}.\footnote{In this paper, we use the terms `generic' and
`typical' interchangeably: the set of coefficient values for which
the anchoring/repelling property fails to hold has measure zero in
the complete coefficient space.} In the same vein, certain dynamics
of spiral waves in anisotropic media, such as phase-locking and
linear drifting of meandering spiral, have been shown to be generic
consequences of rotational symmetry-breaking \cite{R1,R2,LeB}.
\subsubsection*{The basic viewpoint} Consider a piece of cardiac
tissue on which numerous (roughly) circular ablation have been
performed, perhaps in order to treat a patient who is suffering from
atrial fibrillation \cite{Grubb,MDC}.  These surgical procedures
affects both the geometry and the excitability of the tissue. Under
certain modeling assumptions, any system used to model the
electrical activity of the tissue needs to incorporate translational
symmetry-breaking (TSB) components to model the effects of the
circular ablations, and a rotational symmetry-breaking (RSB)
component to model the effects of anisotropy.  Let us model the
electrical properties of such a perturbed piece of anisotropic
cardiac tissue using a modified version of the bidomain equations of
cardiology, under the modeling assumption that the circular ablation
(inhomogeneous) zones consist of a finite number of independent
``sources'' which are localized near distinct sites $\zeta_1,\ldots,
\zeta_n$ in the plane (see \cite{BLM} for a similar hypothesis). The
model then has the form
\begin{align}\label{thebidomain2}
\begin{split} u_t&=\frac{1}{\varsigma}(u-\frac{u^3}{3}-v)+\nabla^2
u+\frac{\alpha \varepsilon}{1+\alpha (1-\varepsilon)}\Psi_{x_1x_1}+
\sum_{j=1}^n\mu_j g_j^{u}(\|x-\zeta_j\|^2,\mu)\\
v_t&=\varsigma(u+\delta-\gamma v)+\sum_{j=1}^n\mu_j g_j^{v}(\|x-\zeta_j\|^2,\mu), \\
\nabla^2& \Psi+\varepsilon g(\alpha,\varepsilon)\Psi_{x_2x_2}=
\varepsilon h(\alpha,\varepsilon) u_{x_2x_2},
\end{split}
\end{align} where $u$ is a transmembrane potential, $v$ controls the recovery of
the action potential, $\Psi$ is an auxiliary potential (without
obvious physical interpretation), $x_1$ is the preferred direction
in physical space in which tissue fibers align, $\varepsilon$~is a
measure of that preference, $g$ and $h$ are appropriate model
functions, $\alpha,\varsigma,\delta$ and $\gamma$ are model
parameters, $\mu=(\mu_1,\ldots,\mu_n)\in \GR^n$ is a small parameter
and $g_j^{u,v}$ are smooth functions, uniformly bounded in their
variables \cite{R2,LR,BEL,Byeah}.\par If the tissue has equal
anisotropy ratios (\textit{i.e.} $\varepsilon=0$) and the
inhomogeneities have no effect on spiral wave dynamics
(\textit{i.e.} $\mu=0$), (\ref{thebidomain2}) decouples into the
FitzHugh-Nagumo equations for $u$ and $v$, and Poisson's equation
for $\Psi$ \cite{R2}.\newl Let $\mathbb{S}\mathbb{E}(2)$ denote the
group of all planar translations and rotations, and fix an integer
$1\leq \jmath^*$ and $\zeta\in \GR^2$. The subgroups
$\mathbb{Z}_{\jmath^*}\sdp \GR^2$ (the notation will be explained
later) and $\GSO(2)_{\zeta}$ of $\GSE(2)$ consist of all cartesian
pairs of translations and rotations about the origin by an integer
multiple of $2\pi/\jmath^*$ radians, and of all rotations about the
point $\zeta$, respectively. Let
$\Gamma=\GC\sdp\mathbb{Z}_{\jmath^*}$ or $\Gamma=\GSO(2)_{\zeta}$.
Then, $\Gamma<\GSE(2)$ and we will say that the semi-flow
$\Psi_{t,\varepsilon,\mu}$ is $\Gamma-$equivariant if it commutes
with the restriction of (\ref{uaction}) to $\Gamma$.\par In the
equivariant dynamical system approach, the particular form of the
functions $g_j^{u,v}$ is unimportant; the analysis is driven by the
fact  that (\ref{thebidomain2}) can sustain spiral wave propagation
\cite{R2,LR,BEL,Byeah} and by the equivariance properties of the
semi-flow $\Phi_{t,\varepsilon,\mu}$ generated by
(\ref{thebidomain2}), namely: if we neglect boundary effects, the
semi-flow \begin{description}
\item[(E1)] is $\GSE(2)-$equivariant when $(\varepsilon,\mu)=0$;
\item[(E2)] is $\GZ_{2}\sdp\GR^2-$equivariant when $\varepsilon\neq 0$ is small and $\mu=0$;
\item[(E3)] preserves rotations around $\zeta_{i}$ (but generically not translations) when
$\varepsilon=0$ and $\mu_{j}=0$ for all $j\neq i$, and
\item[(E4)] is (generically) trivially equivariant when $(\varepsilon,\mu)$ is a generic small parameter vector.
\end{description}
This is but a special case of a more general family of semi-flows
for which \textbf{(E2)} is replaced by the following property: the
semi-flow
\begin{description}
\item[(E2')] is $\GZ_{\jmath^*}\sdp\GR^2-$equivariant when $\varepsilon\neq 0$ is small and
$\mu=0$ for some integer $\jmath^*\geq 1$.
\end{description}
In \cite{BLM,Byeah,Bo1}, we used the dynamical system approach to
study spiral anchoring in media satisfying \textbf{(E1)},
\textbf{(E2')}, \textbf{(E3)} and \textbf{(E4)}: the predictive
power of the method was used to show that in the case $n>1$, spiral
anchoring typically takes place \textit{away} from the
inhomogeneities. At the time, such a statement defied experimental
wisdom.
\subsubsection*{Epicyclic drifting} At this stage, nothing has been said about the nominal topic of this paper:
epicyclic drifting. The various spiral motions observed in
experiments and simulations have been classified according to their
tip path, an arbitrary point on the wave front that is followed in
time \cite{Detal,LOPS}: for instance, the tip path of a (rigidly)
rotating wave is a perfect circle. Barkley \cite{BKT} and Wulff
\cite{Wulff} have shown that the appearance of an epicyclic tip path
can be linked to a 'symmetric Hopf bifurcation:' when that happens,
every spiral wave in the excitable medium is epicyclic.
\par However, other epicyclic behaviour cannot be explained by this
mechanism. When the sizes of the physical domain and of the spiral
core are comparable, the latter is sometimes attracted to the
boundary of the domain and rotates around it in a meandering
fashion. This has been observed in experiments and numerical
simulations in a light-sensitive BZ reaction \cite{YP,ZM}. \par Yet
another instance of epicyclic motion is shown
\begin{figure}[t]
\begin{center}
\includegraphics[width=180pt]{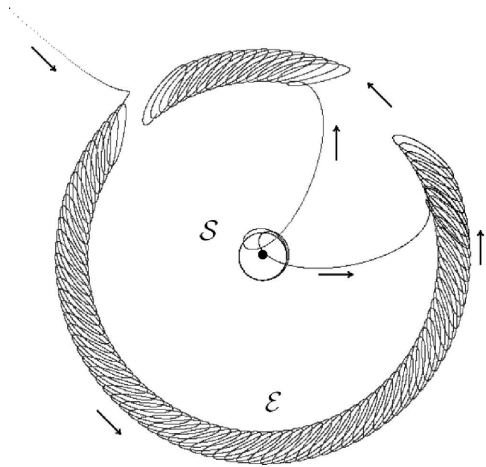}
\caption{Epicyclic motion on the stable epicyclic manifold
$\mathcal{E}$. The arrows indicate the direction of the flow, while
$\mathcal{S}$ corresponds to a repelling (perturbed) rotating wave
solution pinned at the inhomogeneity indicated by the black dot (see
\cite{BLM} for details on spiral anchoring).
\label{epmannn}}\end{center}\hrule
\end{figure}
figure~\thefigure: in a bounded region, all solutions are
attracted/repelled to/by an epicyclic solution manifold. This type
of spiral wave motion is what we refer to as \textit{epicyclic
drifting}. \`A la Poincar\'e-Bendixson, if a system has a stable
epicyclic manifolds (stable in the sense of Lyapunov) it will also
have a repelling rotating wave (see figure~\thefigure), and
\textit{vice-versa}. As such, these manifolds cannot be observed in
fully Euclidean media. What then, can forced Euclidean
symmetry-breaking (FESB) tell us about epicyclic drifting in systems
with the equivariance properties of \textbf{(E1)}, \textbf{(E2')},
\textbf{(E3)} and \textbf{(E4)}? The only work in this vein has been
performed by LeBlanc and Wulff in \cite{LW}, in the case $n=1$ and
without rotational symmetry-breaking: unfortunately, the main tool
in their analysis cannot be used in the general case.
\subsubsection*{Article Overview}
The main object of analysis in the present paper is a
finite-dimensional system of ODE that share the equivariance
properties of \textbf{(E1)}, \textbf{(E2')}, \textbf{(E3)} and
\textbf{(E4)} when $n>1$: it is derived in section \ref{cbe}. Then,
in section~\ref{agat}, we present a preliminary result about
averaging which will subsequently be used to prove our main results:
to wit, when $\varepsilon=0$ or $\jmath^*=1$ and certain conditions
are satisfied, there is a (minimal) parameter wedge region in which
an epicyclic manifold persists. In the case $\jmath^*>1$, the
epicyclic manifold persists in a deleted neighbourhood of the
origin. The parameter wedges are illustrated in \begin{figure}[t]
\begin{center}
\includegraphics[width=360pt]{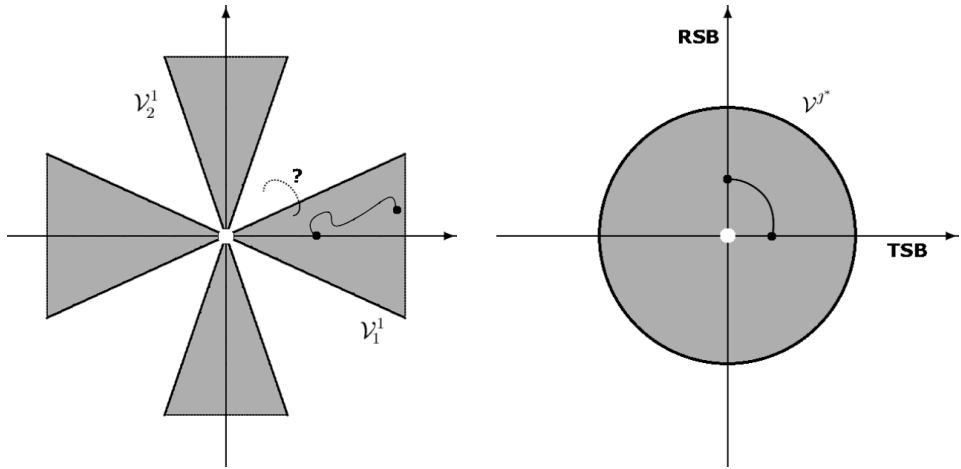}
\caption{On the left, the epicyclic parameter wedge regions for the
case $n=2$ without anisotropy. On the right, the epicyclic deleted
neighbourhood in parameter space for the case $n=1$ with anisotropy
characterized by $\jmath^*>1$. If the semi-flow has an epicyclic
solution manifold for a particular set of parameter values, then the
semi-flow has an epicyclic solution manifold (of the same stability
type) for all parameter values in the adjacent region. Note that
these manifolds continuously deform along any path contained
entirely in the parameter region. The local analysis does not
provide a clear picture of the behaviour as a path leaves a
parameter region. In \cite{Byeah}, for instance, we give an example
where the epicyclic solution manifold disappears as a result of a
saddle-node bifurcation of rotating waves. The notation will be
explained later in the paper. \label{epmannnn}}
\end{center}\hrule
\end{figure} figure~\thefigure, for the case $n=2$. The conditions needed depend on the
kind of forced symmetry-breaking under consideration: in section 3,
we study the general semi-flow under $n$ TSB terms (\textit{i.e.}
$n>1$, $\varepsilon=0$); in section 4, we study the combination of a
single RSB and a single TSB term (\textit{i.e.} $n=1$,
$\varepsilon\neq 0$). We then combine these results in section 5 to
obtain the epicyclic drifting theorems under general FESB. Next, we
perform a simple numerical experiment on \ref{thebidomain2}) with
$n=1$ showing the predicted epicyclic motion: to the best of our
knowledge, the figure in section 6 is the first observed instance of
epicyclic drifting in a numerical simulation of excitable media.
Finally, we give the proofs of two technical results in appendix~A.


\section{Preliminaries} We start with a derivation of the appropriate center bundle equations describing the essential dynamics
of spiral waves near a rotating wave under full Euclidean
symmetry-breaking (FESB). More details on these manipulations can be
found in \cite{Wulff,BLM,SSW1,SSW2,Byeah}. Then, in order to lighten
the text, we introduce some necessary definitions. Finally, we state
an averaging theorem which will be used in later sections of this
work.
\subsection{The Center Bundle
Equations}\label{cbe}  In order to facilitate the subsequent
analysis, we make the same simplifying assumptions and adopt the
same notation as in \cite{BLM,Bo1}.
\par In particular, let $X$ be a Banach space, ${\mathcal
U}\subset\GR\times\mathbb{R}^n$ a neighborhood of the origin, and
$\Phi_{t,\varepsilon,\mu}$ be a smoothly parameterized family
(parameterized by $(\varepsilon,\mu)\in {\mathcal U}$) of smooth
local semi-flows on $X$, and let
\begin{equation}
a:\mathbb{S}\mathbb{E}(2)\longrightarrow \GL(X) \label{a_action}
\end{equation}
be a faithful and isometric representation of
$\mathbb{S}\mathbb{E}(2)$ in the space of bounded, invertible linear
operators on $X$.  For example, if $X$ is a space of functions with
planar domain, a typical $\mathbb{S}\mathbb{E}(2)$ action (such as
(\ref{uaction}) in the preceding section) is given by
$$(a(\gamma)u)(x)=u(\gamma^{-1}(x)),\,\,\,\,\gamma\in\mathbb{S}\mathbb{E}(2).$$
In this paper, we concern ourselves with the study of epicycle
drifting in the case where the two following hypotheses are
satisfied. The first one is simply a re-telling of the equivariance
properties \textbf{(E1)}, \textbf{(E2')}, \textbf{(E3)} and
\textbf{(E4)}, while the second postulates the existence of a
rotating wave in the unperturbed $\GSE(2)-$equivariant semi-flow.
\begin{hyp} \label{hyp1} There exists $1\leq \jmath^*\in \GN$, distinct points $\zeta_1,\ldots,\zeta_n$ in $\GR^2$ such that if $e_j$
denotes the $j^{\mbox{\footnotesize th}}$ vector of the canonical
basis in $\mathbb{R}^n$, then $\forall\,u\in\,X,\ \varepsilon\neq
0,\ \alpha\neq 0,\ t>0,$
\begin{align*}
\Phi_{t,\varepsilon,0}(a(\gamma)u)&=a(\gamma)\Phi_{t,\varepsilon,0}(u)
\iff
\gamma\in\,\mathbb{Z}_{\jmath^*}\sdp \GR^2,\\
\Phi_{t,0,\alpha e_j}(a(\gamma)u)&=a(\gamma)\Phi_{t,0,\alpha
e_j}(u) \iff \gamma\in\,\mathbb{S}\mathbb{O}(2)_{\zeta_j}, \quad\mbox{and}\\
\Phi_{t,0,0}(a(\gamma)u)&=a(\gamma)\Phi_{t,0,0}(u), \quad \forall\,
\gamma\in\,\mathbb{S}\mathbb{E}(2).
\end{align*}
\end{hyp}
\begin{hyp}\label{hyp2}
There exists $u^*\in X$ (with trivial isotropy subgroup) and
$\Omega^*$ in the Lie algebra of $\mathbb{S}\mathbb{E}(2)$ such that
$e^{\Omega^*t}$ is a rotation and
$\Phi_{t,0}(u^*)=a(e^{\Omega^*t})u^*$ for all $t$. Moreover, the set
$\{\,\mu\in\mathbb{C}\,\,|\,\,|\mu|\geq 1\,\}$ is a spectral set for
the linearization $a(e^{-\Omega^*})D\Phi_{1,0}(u^*)$ with projection
$P_*$ such that the generalized eigenspace $\mbox{\rm range}(P_*)$
is three dimensional.
\end{hyp}
As discussed previously, such semi-flows can arise from the family
of perturbed reaction-diffusion systems (1.4) from \cite{BLM} if
$\varepsilon=0$, as well as from the modified bidomain model
(\ref{thebidomain2}) given in the introduction if $\jmath^*=2$, for
instance.
\par It has been shown in \cite{BLM,Bo1,Byeah} that, for small
parameter vectors $(\varepsilon,\mu)\in \GR\times \mathbb{R}^n$, the
essential dynamics of the semi-flow $\Phi_{t,\varepsilon,\mu}$ near
a (hyperbolic) rotating wave is (locally) equivalent to the
semi-flow of the following ordinary differential equations on the
bundle $\mathbb{C}\times\mathbb{S}^1$:
\begin{align}\label{basic12}
\dot{p}&= e^{it}\left[v+\beta G(t,\beta)+\sum_{j=1}^n\lambda_j
H_j((p-\xi_j) e^{-it},\ov{(p-\xi_j)} e ^{it},\lambda)\right]
\end{align}
where $v\in \GC$, $(\beta,\lambda)\in\GR\times \GR^n$,
$\xi_1,\ldots,\xi_n\in \GC$ are all distinct and the functions
$G,H_j$ are smooth, periodic in $t$ and uniformly bounded in $p$,
and $G$ is $2\pi/\jmath^*-$periodic in $t$. The specific form of the
perturbations is a consequence of forced Euclidean symmetry-breaking
from the $\GSE(2)-$equivariance of (\ref{basic12}) under the
following $\GSE(2)-$action on the bundle
$\mathbb{C}\times\mathbb{S}^1$: \begin{align}(x,\theta)\cdot
(p,\varphi)
=(e^{i\theta}p+x,\varphi+\theta),\label{therealaction}\end{align}
for all $(p,\varphi)\in \GC\times \GS^1$ and $(x,\theta)\in
\GSE(2)=\GC\times \GS^1$, where $\GSE(2)=\GC\sdp \GS^1$ with
multiplication $(p_1,\varphi_1) \cdot (p_2,\varphi_2) =
(e^{i\varphi_1}p_2 + p_1, \varphi_1 + \varphi_2).$ The non-standard
multiplication is made explicit by using the semi-direct product
notation $\sdp$. \par Let $\jmath^*\geq 1$ be an integer and
$\xi\in\GC$. In $\GSE(2)=\GC\sdp \GS^1$, the subgroup of rotations
around $\xi$ is given by $$\GS^1_{\xi}=\{(\xi,0)\cdot
(0,\theta)\cdot (-\xi,0):\theta\in \GS^1\},
$$
while the subgroup containing all translations and rotations by
angle $\frac{2\pi k}{\jmath^*}$, $k\in \GZ$, is
$$\GC\sdp\GZ_{\jmath^*}=\left\{\left(x,\frac{2\pi
k}{\jmath^*}\right):k\in \GZ:x\in \GC \right\}.$$
\newpage\noindent Then $\GC\sdp\GZ_{\jmath^*}\simeq \GZ_{\jmath^*}\sdp \GR^2$ and $\GS^1_{\xi}\simeq \GSO(2)_{\zeta}$: under the action
described by (\ref{therealaction}), the center bundle equation
(\ref{basic12}) is
\begin{description}
\item[(C1)] $\GC\sdp\GS^1-$equivariant when $(\beta,\lambda)=0$;
\item[(C2')] $\GC\sdp\GZ_{\jmath^*}-$equivariant when $\beta\neq 0$ is small and $\lambda=0$;
\item[(C3)] $\GS^{1}_{\xi_{\ell}}-$equivariant when
$\beta=0$ and $\lambda_{j}=0$ for all $j\neq \ell$, and
\item[(C4)] (generically) trivially equivariant when $(\beta,\lambda)$ is a generic small parameter vector.
\end{description}
Clearly, (\ref{basic12}) shares the equivariance properties of
Hypothesis 1. As such, $G$ `models' the RSB perturbation while the
various $H_j$ `model' the various TSB terms.
\par It might seem strange that the parameters $(\varepsilon,\mu)$
are replaced by $(\beta,\lambda)$ in (\ref{basic12}), just as the
$\zeta_j\in \GR^2$ are replaced by $\xi_j\in \GC$, but since the
center manifold reduction theorems of \cite{SSW1,SSW2,SSW3,SSW4} do
not provide an explicit relation between the coefficients of the
original system of partial differential equations and the reduced
ordinary differential system of center bundle equations, one cannot
conclude that the parameters are the same in both systems.
\subsection{Definitions} An integral manifold is
\textit{stable} if it has a neighbourhood in which all originating
positive-time solutions approach the manifold exponentially; it is
\textit{hyperbolic} if the linearization of the flow on this
manifold admits no critical eigenvalues. Let $\alpha_0>0$,
$\Delta>0$, $V\subseteq \GR^p$, $\Sigma=\GR\times V\times
[0,\alpha_0]$, $f:\Sigma\to \GR^q$, $g:V\to \GR^q$ and $h:\GR\times
V\to \GR^q$. We say that $f$ is \textit{Lipschitz in Hale's
sense}\label{LipH}, which we denote by $f\in
\Lip(x;\Sigma,\eta(\alpha,V))$, if $f$ is continuous in all of its
arguments and is Lipschitz in $x$ for $(t,x,\alpha)\in \Sigma$ with
continuous Lipschitz constant. \par Next we say that $g$ is
\textit{bounded by $\Delta$ over $V$}, which we denote by $g\in
\mathcal{B}(\Delta;V)$, if $\|g(x)\|\leq \Delta$ for all $x\in V$.
Finally, we denote the fact that $h$ is $T-$periodic in $\phi\in
\GR$ by $h\in\mathfrak{P}^T_{\phi}$. When the sets $\Sigma$ and $V$
are understood from the context, they are omitted. Finally, by abuse
of notation, we shall often denote
$O\left(|x_1|+\cdots+|x_m|\right)$ by $O(x_1,\ldots,x_m)$.
\subsection{A Generalized Averaging Theorem}\label{agat}
Averaging methods are used to determine whether a particular system
has an non-trivial invariant integral manifold by studying an
averaged system. The main theorem is a modified version of one of
Hale's averaging theorems (see \cite{Byeah} for details); it can
easily be extended to the case where $\nu$ is a parameter vector in
$\GR^n$.
\begin{thm}\label{Halesavg}
\textsc{(}$\mathrm{modified\  from\ }$\textsc{\cite{H1}},$\mathrm{\
theorem \ 6.1,\ pp.\ 526-527}$\textsc{)} Let $\sigma_0>0$. Consider
the system of equations
\begin{align}\label{Halesavgsys}
\begin{split}
\dot{x}&=\epsilon \gamma_{\epsilon,\nu} x+\epsilon \Lambda(t,\psi,x,\epsilon,\nu)\\
\dot{\psi}&=d(\epsilon,\nu)+\Theta(t,\psi,x,\epsilon,\nu),
\end{split}
\end{align}
where $\psi\in \GR$, $x\in [-\sigma_0,\infty)$,
$\gamma_{\epsilon,\nu}\neq 0$ depends continuously on
$(\epsilon,\nu)$, and $d$ is defined over
$S_0=[-\epsilon_0,\epsilon_0]\times [-\nu_0, \nu_0]$, with
$d(0,0)=1$. For $\sigma>0$, let
\begin{align*}
\Sigma_{\sigma}&=\GR\times \GR \times [-\sigma,\sigma]\times S_0
\quad\mbox{and}\quad \Sigma_0=\GR\times\GR\times \{0\}\times S_0.
\end{align*}
Suppose $\Theta,\Lambda\in
\mathfrak{P}^{\chi}_t\cap\mathfrak{P}^{\omega}_{\psi}$ and that
\begin{enumerate}[$(i)$]
\item $\Theta$ and $\Lambda$ are real-valued over
$\Sigma_{\sigma_0}$;
\item $\Theta,\Lambda\in \mathcal{B}(\Xi(\epsilon,\nu);\Sigma_0)$ where
  $\Xi(\epsilon,\nu)=O(\epsilon,\nu)$;
\item for all $0\leq \sigma\leq\sigma_0$, $\Theta\in
  \Lip(\psi,x;\Sigma_{\sigma},\theta(\epsilon,\nu,\sigma))$ and
  $\Lambda\in\Lip(\psi,x;\Sigma_{\sigma},\eta(\epsilon,\nu,\sigma))$,
  with
$\theta(\epsilon,\nu,\sigma)=O(\epsilon,\nu,\sigma)$ and
$\eta(\epsilon,\nu,\sigma)=O(\epsilon,\nu,\sigma)$.
\end{enumerate}
Then, there exists $(\epsilon_1,\nu_1)\in (0,\epsilon_0]\times
(0,\nu_0]$ such that for all $$(\epsilon,\nu)\in
S_1=[-\epsilon_1,\epsilon_1]\times [-\nu_1,\nu_1]$$ with
$\epsilon\neq 0$, $(\ref{Halesavgsys})$ has a hyperbolic integral
manifold $\mathcal{T}_{\epsilon,\nu}$ which can be represented as an
invariant torus\index{invariant torus}
$x=\Upsilon_{\epsilon,\nu}(t,\psi)$, where
$$\Upsilon_{\epsilon,\nu}\in \mathcal{B}(D(\epsilon,\nu))\cap
\Lip (\psi,\Omega(\epsilon,\nu))\cap \mathfrak{P}^{\chi}_t\cap
\mathfrak{P}^{\omega}_{\psi},$$ with
$D(\epsilon,\nu),\Omega(\epsilon,\nu)\to 0$ uniformly as
$(\epsilon,\nu)\to 0$. Furthermore, the stability of
$\mathcal{T}_{\epsilon,\nu}$ is exactly determined by the sign of
$\epsilon \gamma_{0,0}$.
\end{thm}

\section{Epicyclic Drifting For $n$ Simultaneous TSB Terms} When $\beta=0$, (\ref{basic12}) gives the dynamics near a hyperbolic
rotating wave for a parameterized family of semi-flows
$\Phi_{t,0,\lambda}$ satisfying the forced-symmetry breaking
conditions in hypothesis \ref{hyp1}. We start with a brief review of
epicyclic drifting in the case $n=1$, which was studied in detail in
\cite{LW}, and then present our new results in the general case
$n>1$.
\subsection{The Case $n=1$} Without loss of generality, we may assume $\xi_1=0$. In this case, the center
bundle equation (\ref{basic12}) reduce to
\begin{align}\label{basiceqs4}
 \dot{p}&={\displaystyle e^{it}\left[v+ \lambda_1
H_1(pe^{-it},\overline{p}e^{it}, \beta)\right]},
\end{align}where $v\in \GC^{\times}$ and $\lambda_1\in\mathbb{R}$ is small. Set
 $\widetilde{H}(w,\overline{w},\lambda_1)=
H_1(w-iv,\overline{w}+i\overline{v},\lambda_1)$.
\begin{thm} \label{nongenericthm} \textsc{(\cite{LW},}$\mathrm{\ re-written\ to\ fit\ the\ current\ symbolism}$\textsc{)} Let
$$I(\rho)=\real\left[\int_0^{2\pi}\!\!e^{-it}\widetilde{H}\left(\rho
    e^{-it},\rho e^{it},0\right)
  dt\right].$$ If $\rho_0>0$ is a hyperbolic solution of $I(\rho)=0$, then for all $\lambda_1\neq 0$ small enough,
  the center bundle equation $(\ref{basiceqs4})$ has an integral (solution) manifold $\mathcal{E}^{1}_{\lambda_1}$ around the origin, whose stability is
  exactly determined by the sign of $\lambda_1 I'(\rho_0)$.
\end{thm}
These solutions represent quasi-periodic motion around the origin in
the $p-$plane and are observable as epicycle-like motion along a
circular boundary in the physical space, with angular frequency
$1+O(\lambda_1)$. Note that the hypotheses of
theorem~\ref{nongenericthm} are not generic: in a random system,
$I(\rho)$ may very well not have a positive hyperbolic root.\par The
presence of a repelling integral manifold could explain the fact
that spirals are sometimes observed to be repulsed by an
inhomogeneity if the spiral tip is located beyond a certain distance
from the perturbation center \cite{MPMPV}.
\subsection{The Case $n>1$} \label{bd} However, the main averaging tool used in \cite{LW}
to obtained theorem (\ref{nongenericthm}) cannot be used to analyze
the situation in the case $n>1$; furthermore, this difficulty yields
a interesting twist, as we shall see in this section.\newl By
re-labeling the indices in (\ref{basic12}) if necessary, we can
temporarily shift our point of view so that $\xi_1$ plays the
central role in the following analysis. Set $\Xi_j=\xi_j-\xi_1$ for
$j=1,\ldots, n$. Then, under the co-rotating frame of reference
$z=p-\xi_1+i e^{i t}v$, (\ref{basic12}) becomes
\begin{align}\label{zdotforced2rescaled}
\dot{z}= e^{i t}\sum_{j=1}^n \lambda_j H_j \big((z-\Xi_j) e^{-i
t}\!\!-i v,\ov{(z-\Xi_j)} e^{i t}\!\!+i\ov{v},\lambda\big).
\end{align}
When $\lambda_1\neq 0$ and $\lambda_2=\cdots=\lambda_n=0$, we find
ourselves in the situation described in the previous subsection.
Now, set $\epsilon=\lambda_1$, $\nu_1=1$ and
$\lambda_j=\nu_j\epsilon$ for $j=2,\ldots, n$, and
$\nu=(\nu_2,\ldots,\nu_n)\in \GR^{n-1}$. Then
(\ref{zdotforced2rescaled}) can be viewed as a perturbation of the
corresponding equation in the case $n=1$. Note that $\Xi_1=0$ and
 $\lambda=(1,\nu)\epsilon$. \newl Equation
(\ref{zdotforced2rescaled}) rewrites as
\begin{equation}\label{system2}
\dot{z}=\epsilon  e^{i t}\sum_{j=1}^n \nu_j H_j \big((z-\Xi_j) e^{-i
t}\!\!-i v,\ov{(z-\Xi_j)} e^{i t}\!\!+i\ov{v},(1,\nu)\epsilon\big).
\end{equation}
Let $\label{funcH} \hat{H}_j(w,\ov{w},\epsilon,\nu)=H_j\big(w-i
v,\ov{w}+i\ov{v},(1,\nu)\epsilon\big)$ for $j=1,\ldots n$. Then
(\ref{system2}) becomes
\begin{align}\label{system3}
\dot{z}=\epsilon  e^{i t}K(z e^{-i t},\ov{z} e^{i t},t,\epsilon,\nu)
\end{align}
where $\displaystyle{K(w,\ov{w},t,\epsilon,\nu)=\sum_{j=1}^n \nu_j
\hat{H}_{j}(w-\Xi_j e^{-i t},\ov{w}-\ov{\Xi}_j e^{i
t},\epsilon,\nu)}$ is $2\pi-$periodic in $t$. Consider the
near-identity change of variables
\begin{align}\label{fcov} z&=w+\epsilon \kappa
(w,\ov{w},t,\epsilon,\nu)
\end{align}
where $\kappa\in \mathfrak{P}^{2\pi}_t$ is differentiable in all of
its variables. Then
\begin{align*} \dot{z}&=\dot{w}+\epsilon \left(\frac{\partial \kappa}{\partial t}+\frac{\partial \kappa}{\partial w }\dot{w}+\frac{\partial \kappa}{\partial \ov{w}}\dot{\ov{w}}\right). 
\end{align*}
Introducing the equivalent complex conjugate equation, this last
system becomes
\begin{equation}\label{laprem}
\left[I_2\index{I2@$I_2$}+\epsilon\begin{pmatrix} \kappa_{w} &
\kappa_{\ov{w}} \\ \ov{\kappa}_w & \ov{\kappa}_{\ov{w}}
\end{pmatrix} \right]\begin{pmatrix} \dot{w} \\
\dot{\ov{w}}\end{pmatrix} =\begin{pmatrix} \dot{z} \\ \dot{\ov{z}}
\end{pmatrix}-\epsilon\begin{pmatrix} \kappa_t \\ \ov{\kappa}_t
\end{pmatrix},
\end{equation}
where $\kappa_{w}$, $\kappa_{\ov{w}}$, $\kappa_t$,
$\ov{\kappa}_{w}$, $\ov{\kappa}_{\ov{w}}$, $\ov{\kappa}_t$ are used
to denote the partial derivatives of $\kappa$ and $\ov{\kappa}$. Set
$$\mathcal{I}=I_2+\epsilon\begin{pmatrix} \kappa_{w} &
\kappa_{\ov{w}} \\ \ov{\kappa}_w & \ov{\kappa}_{\ov{w}}
\end{pmatrix}.$$ Combining (\ref{laprem}) with (\ref{system3})
yields
\begin{equation}\label{theeq}
\begin{pmatrix} \dot{w} \\ \dot{\ov{w}}\end{pmatrix}=
\epsilon\mathcal{I}^{-1}
\begin{pmatrix} e^{it}K\left((w+\epsilon\kappa)e^{-it},(\ov{w}+\epsilon\ov{\kappa})e^{it},t,\epsilon,\nu\right)- \kappa_t \\ e^{-it}\ov{K}\left((w+\epsilon\kappa)e^{-it},(\ov{w}+\epsilon\ov{\kappa})e^{it},t,\epsilon,\nu\right)- \ov{\kappa}_t \end{pmatrix}
\end{equation}
By Taylor's theorem, there are appropriate continuous bounded
functions $A_1,$ $A_2$ and $A_3\in \mathfrak{P}^{2\pi}_t$ satisfying
\begin{align*}
e^{it}K\left((w+\epsilon\kappa)e^{-it},(\ov{w}+\epsilon\ov{\kappa})e^{it},t,\epsilon,\nu\right)&=e^{it}K\left(we^{-it},\ov{w}e^{it},t,0,\nu\right)
+\epsilon
A_1(w,\ov{w},t,\epsilon,\nu) \\
\kappa_t(w,\ov{w},t,\epsilon,\nu)&=\kappa_t(w,\ov{w},t,0,\nu)+\epsilon
A_2(w,\ov{w},t,\epsilon,\nu)
\end{align*}
and
$$\mathcal{I}^{-1}=
\begin{pmatrix}
1-\epsilon \kappa_{w}^0 & -\epsilon \kappa_{\ov{w}}^0 \\
-\epsilon \ov{\kappa}_{w}^0 & 1-\epsilon \ov{\kappa}_{\ov{w}}^0 \\
\end{pmatrix} +\epsilon^2A_3(w,\ov{w},t,\epsilon,\nu),$$ where
$$\kappa_{w}^0=\kappa_w(w,\ov{w},t,0,\nu)\quad \mbox{and} \quad
\kappa_{\ov{w}}^0=\kappa_{\ov{w}}(w,\ov{w},t,0,\nu).$$\normalsize
With these, (\ref{theeq}) re-writes (upon dropping the equivalent
complex conjugate equation) as
\begin{align}
\label{listen} \dot{w}&=\epsilon
\left(e^{it}K(we^{-it},\ov{w}e^{it},t,0,\nu)-\kappa_t(w,\ov{w},t,0,\nu)\right)+\epsilon^2
\mathcal{H}(w,\ov{w},t,\epsilon,\nu),
\end{align}
where $\mathcal{H}\in \mathfrak{P}^{2\pi}_t$ is bounded and
continuous in all its variables. Denote the average value of
\label{pappequiv}$e^{it}K(we^{-it},\ov{w}e^{it},t,0,\nu)$ by
\begin{equation}\label{avgM}M^1(w,\ov{w},\nu)=\frac{1}{2\pi}\int_{0}^{2\pi} \!\!\!\!e^{it}K(we^{-it},\ov{w}e^{it},t,0,\nu)\, dt.\index{average value!Ma1@$M^1$}\index{Maaaa1@$M^1$}\end{equation} Then $$e^{it}K(we^{-it},\ov{w}e^{it},t,0,\nu)= M^1(w,\ov{w},\nu)+F(w,\ov{w},t,\nu),$$ where $F\in \mathfrak{P}^{2\pi}_t$ is uniformly continuous and \begin{equation}\label{noconstant}\int_{0}^{2\pi}\!\!\!\!F(w,\ov{w},t,\nu)\, dt=0.\end{equation} Let $\kappa$ be an antiderivative of $F$ with respect to $t$. Then $\kappa\in \mathfrak{P}^{2\pi}_t$ by (\ref{noconstant}) and $$F(w,\ov{w},t,\nu)-\kappa_t(w,\ov{w},t,0,\nu)=0.$$ With such a $\kappa$, (\ref{listen}) simplifies to
\begin{align}
\label{avg} \dot{w}&=\epsilon M^1(w,\overline{w},\nu)+\epsilon^2
\mathcal{H}(w,\ov{w},t,\epsilon,\nu).
\end{align}
It is easy to see that $M^1(w,\overline{w},0)$ is
$\GS^1-$equivariant (see appendix for details); as such, there is a
continuous function $L_1:\GR\to \GC$ such that
$M^1(w,\ov{w},0)=wL_1(w\ov{w})$ \cite[p.~360]{GSS}.
\par By Taylor's theorem, there are appropriate continuous bounded
functions $M_j$, for $j=2,\ldots, n$, such that
$$M^1(w,\ov{w},\nu)=M^1(w,\ov{w},0)+\sum_{j=2}^n\nu_j
M_{j}(w,\ov{w},\nu)$$ and so (\ref{avg}) becomes
\begin{equation}\label{avg2}
\dot{w}=\epsilon wL_1(w\ov{w})+\epsilon W(w,\ov{w},t,\epsilon,\nu),
\end{equation}
where
\begin{equation}\label{funcW}W(w,\ov{w},t,\epsilon,\nu)=\sum_{j=2}^n\nu_j
M_j(w,\ov{w},\nu)+\epsilon
\mathcal{H}(w,\ov{w},t,\epsilon,\nu).\end{equation}
Differentiating the polar coordinates $w=\rho e^{-i(\psi-t)}$ yields 
\begin{align*}
\dot{\rho}&= \real\left[\dot{w}e^{i(\psi-t)} \right]\\
\dot{\psi}&= 1-\frac{1}{\rho}\imag\left[\dot{w}e^{i(\psi-t)}\right].
\end{align*}
But \begin{align*}
\dot{w}e^{i(\psi-t)}&=\left(\epsilon wL_1(w\ov{w})+\epsilon W(w,\ov{w},t,\epsilon,\nu)\right)e^{i(\psi-t)}\\
&=\left(\epsilon \rho e^{-i(\psi-t)}L_1(\rho^2)+\epsilon W(\rho e^{-i(\psi-t)},\rho e^{i(\psi-t)},t,\epsilon,\nu)\right)e^{i(\psi-t)}\\
&=\epsilon\rho L_1(\rho^2)+\epsilon e^{i(\psi-t)} W(\rho
e^{-i(\psi-t)},\rho e^{i(\psi-t)},t,\epsilon,\nu)
\end{align*}
and so
\begin{align}\label{polavg}
\begin{split}
\dot{\rho}&= \epsilon R^1_0(\rho)+ \epsilon R(t,\psi,\rho,\epsilon,\nu)\\
\dot{\psi}&= 1+\epsilon\Psi_0(\rho)+\epsilon
\Psi(t,\psi,\rho,\epsilon,\nu),
\end{split}
\end{align}
where $R^1_0(\rho)= \rho\real \left[L_1(\rho^2)\right]$,
$\Psi_0(\rho)=-\imag\left[L_1(\rho^2)\right]$\index{Ra10@$R^1_0$}\index{drifting!function}\index{boundary
drifting!function} and
\begin{align}
\label{thewhat2}
\begin{split}
R(t,\psi,\rho,\epsilon,\nu)&=\real\left[e^{i(\psi-t)} W(\rho e^{-i(\psi-t)},\rho e^{i(\psi-t)},t,\epsilon,\nu)\right]\\
\Psi(t,\psi,\rho,\epsilon,\nu)&=-\frac{1}{\rho}\imag\left[e^{i(\psi-t)}
W(\rho e^{-i(\psi-t)},\rho e^{i(\psi-t)},t,\epsilon,\nu)\right].
\end{split}
\end{align}
Note that $R,\Psi\in
\mathfrak{P}^{2\pi}_t\cap\mathfrak{P}^{2\pi}_{\psi}$ and that $\Psi$
is not defined at $\rho=0$.  We now give sufficient conditions for
the existence of an integral manifold in (\ref{polavg}).
\begin{thm}\label{thmtorus} Assume that $R$ and $\Psi$, as defined in $(\ref{thewhat2})$, are $C^1$ on
intervals away from $\rho=0$ and that the averaged equation
\begin{equation}\label{averaged}
\dot{\rho}=\epsilon R^1_0(\rho)
\end{equation}
has an equilibrium $\rho_1>0$ with
$D_{\rho}R_0^1(\rho_1)=\gamma_1\neq 0$. If the parameters are small
enough to satisfy the conditions outlined in the proof below, then
 $(\ref{polavg})$ has an invariant torus $\hat{\cal
E}_{\epsilon,\nu}$, whose stability is exactly determined by the
sign of $\epsilon\gamma_1$.
\end{thm}
  \proof By the implicit function theorem, there is a neighbourhood
$$U=(-\epsilon_*,\epsilon_*)\times \prod_{j=2}^n (-\nu_{j,*},\nu_{j,*})$$
in parameter space and a continuous function $\rho:U\to \GR^+$ such
that $\rho(0,0)=\rho_1$,
$$\epsilon R_0^1(\rho(\epsilon,\nu))\equiv
0\quad\mbox{and}\quad
D_{\rho}R_0^1(\rho(\epsilon,\nu))=\gamma_{\epsilon,\nu}\neq 0,$$
where $\gamma_{\epsilon,\nu}\gamma_1>0$ for all $(\epsilon,\nu)\in
U$, \textsl{i.e.} the stability of the equilibria
$\rho(\epsilon,\nu)$ is the same as that of $\rho_0$ for all
$(\epsilon,\nu)\in U$. \par When $\epsilon=0$, the phase space of
(\ref{polavg}) is foliated by invariant tori and so, from now on, we
will assume that $\epsilon\neq 0$. Consider the change of variables
$\rho=\rho(\epsilon,\nu)+x$ in (\ref{polavg}). Differentiating the
new coordinates, we get $\dot{x}=\dot{\rho}$ and the equivalent
system
\begin{align*}
\begin{split}
\dot{x}&=\epsilon R_0^1(\rho(\epsilon,\nu)+x)+\epsilon R(t,\psi,\rho(\epsilon,\nu)+x,\epsilon,\nu)\\
\dot{\psi}&=1+\epsilon\Psi_0(\rho(\epsilon,\nu)+x)+\epsilon\Psi(t,\psi,\rho(\epsilon,\nu)+x,\epsilon,\nu).
\end{split}
\end{align*}
By Taylor's theorem, there are continuously differentiable functions
$B_1$ and $B_2$ such that
\begin{align*}
R_0^1(\rho(\epsilon,\nu)+x)&=R_0^1(\rho(\epsilon,\nu))+D_{x} R_0^1(\rho(\epsilon,\nu))x+B_1(x,\epsilon,\nu)x^2\\
\Psi_0(\rho(\epsilon,\nu)+x)&=\Psi_0(\rho(\epsilon,\nu))+
B_2(x,\epsilon,\nu)x.
\end{align*}
Since $R_0^1(\rho(\epsilon,\nu))\equiv 0$ and $D_{x}
R_0^1(\rho(\epsilon,\nu))=\gamma_{\epsilon,\nu}$, we obtain the new
system
\begin{align}
\begin{split}\label{thexeqs}
\dot{x}&=\epsilon\gamma_{\epsilon,\nu}x+\epsilon\Lambda(t,\psi,x,\epsilon,\nu) \\
\dot{\psi}&=d(\epsilon,\nu)+\Theta(t,\psi,x,\epsilon,\nu),
\end{split}
\end{align}
where
\begin{align*}
\Lambda(t,\psi,x,\epsilon,\nu)&=B_1(x,\epsilon,\nu)x^2+R(t,\psi,\rho(\epsilon,\nu)+x,\epsilon,\nu)\\
\Theta(t,\psi,x,\epsilon,\nu)&=\epsilon B_2(x,\epsilon,\nu)x+\epsilon\Psi(t,\psi,\rho(\epsilon,\nu)+x,\epsilon,\nu)\\
d(\epsilon,\nu)&=1+\epsilon\Psi_0(\rho(\epsilon,\nu))
\end{align*}
are at least $C^1$ by hypothesis.\par Let $U^+=\{\varsigma\in U:
\varsigma_i>0 \mbox{ for all }i=1,\ldots,n\}$ and
$(\epsilon_0,\nu_0)\in U^+$. Define
$$S_0=[-\epsilon_0,\epsilon_0]\times
\prod_{j=2}^n[-\nu_{0,j},\nu_{0,j}].$$  As $\rho_1>0$ and
$\rho(\epsilon,\nu)$ is continuous on $U$, it is possible to chose
$(\epsilon_0,\nu_0)$ in such a way that
 $$ \sigma_0=\min_{(\epsilon,\nu)\in S_0}\left\{\rho(\epsilon,\nu)\right\}-\textstyle{\frac{1}{2}}\rho_1>0. $$  If $x\geq -\sigma_0$, then
$\rho=\rho(\epsilon,\nu)+x\geq \rho(\epsilon,\nu)-\sigma_0\geq
\textstyle{\frac{1}{2}}\rho_1$ for all $(\epsilon,\nu)\in S_0$. In
that case, $\Theta$ and $\Lambda$ are continuously differentiable,
as $R$ and $\Psi$ are continuously differentiable in $\rho$ on
$[\frac{1}{2}\rho_1,\infty)$. Note further that $\Theta,\Lambda\in
\mathfrak{P}^{2\pi}_t\cap \mathfrak{P}^{2\pi}_{\psi}$. \par Set
$\Sigma_0=\GR\times \GR\times \{0\}\times S_0$, and
$\Sigma_{\sigma}=\GR\times\GR\times [-\sigma,\sigma]\times S_0.$
Then \label{pappH1H2}
\begin{enumerate}
\item $\Theta$ and $\Lambda$ are bounded by a function $\Xi(\epsilon,\nu)=O(\epsilon,\nu_2,\ldots,\nu_n)$ over $\Sigma_0$ (see appendix for details), and
\item for all $0\leq \sigma\leq \sigma_0$, $\Theta$ and $\Lambda$ are Lipschitz in Hale's sense (with Lipschitz con\-stants $\theta(\epsilon,\nu,\delta)=O(\epsilon,\nu_2,\ldots,\nu_n,\delta)$ and $\eta(\epsilon,\nu,\delta)=O(\epsilon,\nu_2,\ldots,\nu_n,\delta)$, respectively) over~$\Sigma_{\sigma}$ (see appendix for details).\end{enumerate}
Accordingly, theorem \ref{Halesavg} can be applied to show there is
a neighbourhood $S_1\subseteq S_0$ of the origin in parameter space
for which (\ref{polavg}) (since it is equivalent to (\ref{thexeqs}))
has an invariant torus $\hat{\mathcal{T}}_{\epsilon,\nu}$ when
$(\epsilon,\nu)\in S_1$. Furthermore, the stability of
$\hat{\mathcal{T}}_{\epsilon,\nu}$ is the same as that of the
hyperbolic equilibrium $\rho(\epsilon,\nu)$, which is given by
$\epsilon\gamma_1$. \qed\newl \normalsize The invariant torus
$\hat{\mathcal{T}}_{\epsilon,\nu}$ appearing in the proof of
Theorem~\ref{thmtorus} can be parameterized by a relation of the
form $x=\Upsilon_{\epsilon,\nu}(\theta_1,\theta_2)$, where
$\theta_1,\theta_2\in \GS^1$. Let
\begin{align} \big\langle\hat{\mathcal{T}}_{\epsilon,\nu}\big\rangle=\frac{1}{4\pi^2}\int_{0}^{2\pi}\!\!\!\!\int_{0}^{2\pi}\!\!\!\! \Upsilon_{\epsilon,\nu}(\theta_1,\theta_2)\,d\theta_1d\theta_2\end{align} denote the \textit{center} of $\hat{\cal
T}_{\epsilon,\nu}$, and let $\hat{\cal E}_{\epsilon,\nu}$ be the
corresponding \textit{epicyclic manifold} of (\ref{system3}), in
which all solutions are epicycles when projected upon the $z-$plane.
\par Define the average value
\begin{multline}
\index{[]D@$[-]_{\DD}$}[\hat{\cal
E}_{\epsilon,\nu}]_{\DD}=\frac{1}{4\pi^2}\int_{0}^{2\pi}\!\!\!\!\int_{0}^{2\pi}\!\!\!\!\big(\big(\rho(\epsilon,\nu)+\big\langle\hat{\mathcal{T}}_{\epsilon,\nu}\big\rangle\big)e^{-i(\psi-t)}
\\ +\epsilon\kappa\big(\big(\rho(\epsilon,\nu)+\big\langle\hat{\mathcal{T}}_{\epsilon,\nu}\big\rangle\big)e^{-i(\psi-t)},\mbox{c.c.},t,\epsilon,\nu\big) \big)\, d\psi dt.\label{cd}\end{multline} If $\mathcal{\hat{T}}_{\epsilon,\nu}$ is stable (in the sense of theorem \ref{Halesavg}), we shall say that $[\hat{\mathcal{E}}_{\epsilon,\nu}]_{\DD}$ is the \textit{center of drifting}\label{cdriftoo} of $\hat{\mathcal{E}}_{\epsilon,\nu}$.
 \begin{thm}\label{thmnld1} Suppose the hypotheses of theorem~$\ref{thmtorus}$ are satisfied. Then there exists a wedge-shaped region near
    $\lambda=0$ of the form
 $$
{\cal V}_{1}=\{(\lambda_1,\ldots,\lambda_n)\in
\GR^n\,:\,|\lambda_j|<V_{1,j}|\lambda_1|,\,\,\,V_{1,j}>0,\,\,\mbox{\rm
for $j\neq 1$ and $\lambda_1$ near}\,\,0\,\}
 $$
such that for all $0\neq \lambda\in {\cal V}_{1}$, $(\ref{basic12})$
has an epicycle manifold ${\cal E}_{\lambda}^1$, with $[{\cal
E}^1_{\lambda}]_{\DD}$ near, but generically not at, $\xi_1$.
Furthermore, $[{\cal E}^1_{\lambda}]_{\DD}$ is a center of drifting
when $\lambda_1\gamma_1<0$.
\end{thm}
  \proof According to theorem~\ref{thmtorus}, there are constants
$\epsilon_1,\nu_{1,2},\ldots,\nu_{1,n}>0$ and a neighbourhood
$$S_1=[-\epsilon_1,\epsilon_1]\times
\prod_{j=2}^n[-\nu_{1,j},\nu_{1,j}]$$ such that (\ref{polavg}) has
an integral manifold $\hat{\mathcal{E}}_{\epsilon,\nu}$ whenever
$(\epsilon,\nu)\in S_1$. \par\noindent For $j\neq 1$, set
$\lambda_1=\epsilon\neq 0$, $\lambda_j=\nu_j\epsilon$ and
$V_{1,j}=\nu_{1,j}$. Then $\lambda\in {\cal V}_1$ as
$$|\lambda_j|\leq |\nu_j|\cdot |\lambda_1|\leq
V_{1,j}|\lambda_1|\quad \mbox{for }j\neq 1,$$ and (\ref{basic12})
has an integral manifold ${\cal
E}_{\lambda}^1=\xi_1-ie^{it}v+\hat{\mathcal{E}}_{\epsilon,\nu}$.
Furthermore, $[{\cal E}_{\lambda}^1]_{\DD}=\xi_1+O(\lambda_1)$ and
so $[\mathcal{E}_{\lambda}^1]_{\DD}\neq \xi_1$ for a generic
$0\neq\lambda\in \mathcal{V}_1$. The conclusion on the stability of
${\cal E}^1_{\lambda}$ then follows directly from
theorem~\ref{thmtorus}. \qed\normalsize \begin{rem}
\begin{enumerate}\item These
isolated epicycle manifolds need not in general be unique for a
given $\lambda\in {\cal V}_1$ as $R_0^1(\rho)=0$ may have any number
of hyperbolic solutions. \item In generic semi-flows, all that can
be said with certainty from the analysis when the parameter values
stray outside of ${\cal V}_{1}$ is that the epicycle manifolds in
(\ref{basic12}) drift away from $\xi_1$, which cannot then be a
center of drifting. This is not unlike the situation with regards to
spiral anchoring \cite{BLM}. Richer dynamics and interactions with
rotating waves can also take place; for instance in \cite{Byeah}, we
gave an example in which the epicyclic manifold collapses at a
saddle-node bifurcation of rotating waves. \item Note that the
actual parameter region in which epicyclic drifting is observed may
be much larger than $\mathcal{V}_1$: however, our local analysis
cannot be used to obtain global results.
\end{enumerate}\end{rem} The preceding results have been achieved by
considering (\ref{basic12}) under a co-rotating frame of reference
around $\xi_1$. Of course, since the choice for $\xi_1$ was
arbitrary, corresponding results must also be achieved, in exactly
the same manner, when the viewpoint shifts to another~$\xi_k$.
Indeed, for $j=1,\ldots,n$, define the average functions
$$M^{j}(w,\ov{w})=\frac{1}{2\pi}\int_{0}^{2\pi}\!\!\!\!e^{it}\hat{H}_j(we^{-it},\ov{w}e^{it},0,0)\,
dt;$$ as before, each $M^j$ is $\GS^1-$equivariant and so there are
continuous functions $L_j:\GR\to \GC$ such that
$M^{j}(w,\ov{w})=wL_j(w\ov{w})$. We will call
$$R^{j}_0(\rho)=\rho\real \left[L_j(\rho^2)\right]$$ the
\textit{epicycle functions of} (\ref{basic12}).
\begin{cor}\label{thmnld2} Let $k\in \{1,\ldots, n\}$. If $\rho_*>0$ is such that $$R_0^k(\rho_*)=0 \quad\mbox{and}\quad D_{\rho}R_0^k(\rho_*)=\gamma_*\neq 0,$$
then there exists a wedge-shaped region near
    $\lambda=0$ of the form
 $$
{\cal V}_{k}=\{(\lambda_1,\ldots,\lambda_n)\in
\GR^n\,:\,|\lambda_j|<V_{k,j}|\lambda_k|,\,\,\,V_{k,j}>0,\,\,\mbox{\rm
for $j\neq k$ and $\lambda_k$ near}\,\,0\,\}
 $$
such that for all $0\neq \lambda\in {\cal V}_{k}$, $(\ref{basic12})$
has an epicycle manifold ${\cal E}_{\lambda}^k$, with $[{\cal
E}^k_{\lambda}]_{\DD}$ near, but generically not at, $\xi_k$.
Furthermore, $[{\cal E}^k_{\lambda}]_{\DD}$ is a center of drifting
when $\lambda_k\gamma_*<0$.
\end{cor}
  \proof The epicycle function $R_0^k$ is exactly the function that would appear in
(\ref{averaged}) had the preceding analysis been done around
$\xi_k$. Theorems~\ref{thmtorus} and~\ref{thmnld1}\ can then be
applied directly to obtain the desired result.\qed\newl Clearly, the
remarks appearing after the proof of theorem~\ref{thmnld1} still
hold. There is one last statement to be made concerning epicycle
manifolds: theorem~\ref{thmnld2} only gives sufficient conditions
for their existence in (\ref{basic12}). In section~\cite{Byeah}, we
have provided an example that shows that they are not, in fact,
necessary conditions.


\section{Epicyclic Drifting For Combined RSB-TSB Terms} In this section, we investigate another way in which the Euclidean
symmetry can be broken: by combining rotational and translational
symmetry breaking. In effect, we are lifting the restriction
$\beta=0$, with $n=1$ in (\ref{basic12}). \par It turns out that the
value of $\jmath^*$ plays a crucial role in the analysis: the cases
$\jmath^*=1$ and $\jmath^*>1$ are essentially different. The general
lines are very similar to those of the preceding section, as such,
the proofs are omitted in order to avoid tedious repetitions. In
either case, however, me assume without loss of generality that
$\xi_1=0$.
\subsection{The Case $\jmath^*=1$} Let $F_G:\GR\times\GR^2\to \GC$ be defined by
$$ F_G(t,\beta)=e^{it}\Big[-iv+\beta \sum_{m\neq
-1}\,\frac{g_m(\beta)e^{im t}}{i(m+1)}\Big],$$ where the
$g_{m}(\beta)$ are the Fourier coefficients of $G\in
\mathfrak{P}^{2\pi}_t$. Set $z=p-F_G$. Then, (\ref{basic12})
rewrites as
\begin{align}
\label{zdotforcedyeah1}
\begin{split}
\dot{z}=\beta g_{-1}(\beta)+\beta e^{it}
H((z+F_G(t,\beta)e^{-it},\mbox{c.c.},\lambda_1),
\end{split}
\end{align}
where \mbox{c.c.} represents throughout the complex conjugate of the
preceding term. Generically, $g_{-1}(0)\neq 0$. Set
$\epsilon=\beta$, $\nu=\lambda_1$ and
$\epsilon=\hat{\epsilon}\lambda_1$. Then (\ref{zdotforcedyeah1})
transforms to
\begin{equation}
\dot{z}=\nu H_*(ze^{-it},\ov{z}e^{it},t,\hat{\epsilon},\nu),
\label{newzdotforced1}
\end{equation}
where
\begin{equation}\label{thefirstH}H_*(w,\ov{w},t,\hat{\epsilon},\nu)=\hat{\epsilon}g_{-1}(\hat{\epsilon}\nu)+e^{it}
H(w+F_G(t,\hat{\epsilon}\nu)e^{-it},\ccc,\nu)\end{equation} is
$2\pi-$periodic in $t$, smooth and uniformly bounded in $w$.
Consider the near-identity change of variables
\begin{align}\label{ahahfcov}
z&=w+\nu \varrho (w,\ov{w},t,\epsilon,\nu)
\end{align}
where $\varrho\in \mathfrak{P}^{2\pi}_t$ is continuous in all of its
variables and to be determined later. This change of variables
transform (\ref{newzdotforced1}) into the equivalent system
\begin{align}
\label{ahahlisten} \dot{w}&=\nu
\left(H_*(we^{-it},\ov{w}e^{it},t,\hat{\epsilon},0)-\varrho_t(w,\ov{w},t,\hat{\epsilon},0)\right)+\nu^2
\mathcal{H}_*(w,\ov{w},t,\hat{\epsilon},\nu),
\end{align}
where $\mathcal{H}_*\in \mathfrak{P}^{2\pi}_t$ is bounded and
continuous in all its variables. Denote the average value of
\label{ahahpappequiv}$H_*(we^{-it},\ov{w}e^{it},t,\hat{\epsilon},0)$
by
\begin{equation*}\label{ahahavgM}M_*(w,\ov{w},\hat{\epsilon})=\frac{1}{2\pi}\int_{0}^{2\pi} \!\!\!\!H_*(we^{-it},\ov{w}e^{it},t,\hat{\epsilon},0)\, dt.\end{equation*}
Since
\begin{align*}
H_*(we^{-it},\ov{w}e^{it},t,\hat{\epsilon},0)&=\hat{\epsilon}g_{-1}(0)+e^{it}H\left((w+F_G(t,0))e^{-it}, \ccc,0\right) \\
&=\hat{\epsilon}g_{-1}(0)+e^{it}H(we^{-it}\!\!-iv,\ov{w}e^{it}\!\!+i\ov{v},0),
\end{align*} we have $$M_*(w,\ov{w},\hat{\epsilon})=\hat{\epsilon}g_{-1}(0)+\frac{1}{2\pi}\int_{0}^{2\pi}\!\!\!\!e^{it} H(we^{-it}\!\!-iv,\ov{w}e^{it}\!\!+i\ov{v},0)\, dt.$$ Then $$H_*(we^{-it},\ov{w}e^{it},t,\hat{\epsilon},0)= M_*(w,\ov{w},\hat{\epsilon})+F_*(w,\ov{w},t,\hat{\epsilon}),$$ where $F_*\in \mathfrak{P}^{2\pi}_t$ is uniformly continuous and \begin{equation}\label{ahahnoconstant}\int_{0}^{2\pi}\!\!\!\!F(w,\ov{w},t,\hat{\epsilon})\, dt=0.\end{equation} Let $\varrho$ be an antiderivative of $F_*$ with respect to $t$. Then $\varrho\in \mathfrak{P}^{2\pi}_t$ by (\ref{ahahnoconstant}) and $$F_*(w,\ov{w},t,\hat{\epsilon})-\varrho_t(w,\ov{w},t,\hat{\epsilon},0)=0.$$ With such a $\varrho$, (\ref{ahahlisten}) simplifies to
\begin{align}
\label{ahahavg} \dot{w}&=\nu
M_*(w,\overline{w},\hat{\epsilon})+\nu^2
\mathcal{H}_*(w,\ov{w},t,\hat{\epsilon},\nu).
\end{align}
As $M_*(w,\overline{w},0)$ is also $\GS^1-$equivariant, there is a
continuous function $L_*:\GR\to \GC$ such that
$M_*(w,\ov{w},0)=wL_*(w\ov{w})$, and so (\ref{ahahavg}) becomes
\begin{equation}\label{ahahavg2}
\dot{w}=\nu wL_*(w\ov{w})+\nu W_*(w,\ov{w},t,\hat{\epsilon},\nu),
\end{equation}
where
\begin{equation}\label{ahahfuncW}W_*(w,\ov{w},t,\hat{\epsilon},\nu)=\hat{\epsilon}g_{-1}(0)+\nu
\mathcal{H}_*(w,\ov{w},t,\hat{\epsilon},\nu).\end{equation}
Differentiating the polar coordinates $w=\rho e^{-i(\psi-t)}$ and
substituting in \ref{ahahavg2} yields
\begin{align}\label{ahahpolavg}
\begin{split}
\dot{\rho}&= \nu R_*(\rho)+ \nu R(t,\psi,\rho,\hat{\epsilon},\nu)\\
\dot{\psi}&= 1+\nu\Psi_*(\rho)+\nu
\Psi(t,\psi,\rho,\hat{\epsilon},\nu),
\end{split}
\end{align}
where $R_*(\rho)= \rho\real \left[L_*(\rho^2)\right]$,
$\Psi_*(\rho)=-\imag\left[L_*(\rho^2)\right]$,\index{Ra*@$R_*$}\index{drifting!function}\index{boundary
drifting!function} $R,\Psi\in \mathfrak{P}^{2\pi}_t\cap
\mathfrak{P}^{2\pi}_{\psi}$ and $\Psi$ is not continuous at
$\rho=0$. We now provide sufficient conditions for the existence of
an epicycle manifold in (\ref{basic12}).
\begin{thm}\label{ahahthmtorus} Assume that $R$ and $\Psi$, as defined in $(\ref{ahahpolavg})$, are $C^1$ on
intervals away from $\rho=0$ and that the averaged equation
\begin{equation}\label{ahahaveraged}
\dot{\rho}=\epsilon R_*(\rho)
\end{equation}
has an equilibrium $\rho_0>0$ with $D_{\rho}R_*(\rho_0)=\gamma_0\neq
0$. If the parameters are small enough, there exists a wedge-shaped
region near
    $(\beta,\lambda_1)=(0,0)$ of the form
$${\cal V}=\{(\beta,\lambda_1)\in \GR^2\,:\,|\beta|<K|\lambda_1|,\,\,\,K>0,\,\,\mbox{$\lambda_1$ near}\,\,0\,\}
$$
such that for all $(\beta,\lambda_1)\in {\cal V}$, $\beta\neq 0$,
$(\ref{basic12})$ has an epicycle manifold ${\cal
G}_{\beta,\lambda_1}^1$, with $[{\cal G}^1_{\beta,\lambda_1}]_{\DD}$
near, but generically not at, the origin. Furthermore, $[{\cal
G}^1_{\beta,\lambda_1}]_{\DD}$ is a center of drifting when
$\lambda_1\gamma_0<0$.
\end{thm}
The remarks after theorem \ref{thmnld1} still hold after having been
suitably modified.
\subsection{The Case $\jmath^*>1$} The case $\jmath^*>1$ is handled
slightly differently.  Let $C^{0}_{\GR}(\GC)$ and $C^{1}_{\GR}(\GC)$
be the spaces of continuous and continuously differentiable
functions from $\GR$ to $\GC$, respectively. Then
$$C^0_{2\pi/{\jmath^*}}=\{f: f\in \mathfrak{P}^{2\pi/{\jmath^*}}_t\cap
C^0_{\GR}(\GC)\}\quad\mbox{and}\quad C^1_{2\pi/{\jmath^*}}=\{f:f\in
\mathfrak{P}^{2\pi/{\jmath^*}}_t\cap C^1_{\GR}(\GC)\}$$ are Banach
spaces when endowed with the respective norms
$$||u||_0=\sup\{|u(x)|:x\in [0,\textstyle{2\pi/{\jmath^*}}]\}\quad
\mbox{and}\quad ||u||_1=||u||_0+||u'||_0,$$ and the linear operator
${\cal Y}:C^1_{2\pi/{\jmath^*}}\to C^0_{2\pi/{\jmath^*}}$ defined by
${\cal Y}(u)=iu+u'$ is bounded, invertible and has bounded inverse.
\newl Define the nonlinear operator ${\cal
H}_G:C^1_{2\pi/{\jmath^*}}\times {\mathbb R}^2\to
C^0_{2\pi/{\jmath^*}}$ by
\begin{equation}\label{HG}
{\cal H}_G(u,\beta,\lambda_1)={\cal Y}(u)-\lambda_1 H\left(u-iv+
\beta\sum_{m\in {\mathbb Z}}\,\frac{g_m(\beta)e^{im{\jmath^*}
t}}{i(m{\jmath^*}+1)}, \mbox{c.c.},\lambda_1\right),
\end{equation} where the $g_{m}(\beta)$ are as they were in the previous section. But ${\cal H}_G(0,0,0)=0$
and $D_1{\cal H}_G(0,0,0)=i\neq 0$ and so, by the implicit function
theorem, there is a neighbourhood ${\cal N}$ of the origin in
${\mathbb R}^2$ and a unique smooth function $U:{\mathbb R}^2\to
C^1_{2\pi/{\jmath^*}}$ satisfying $U(0,0)=0$ and ${\cal
H}_G(U(\beta,\lambda_1),\beta,\lambda_1)\equiv 0$ for all
$(\beta,\lambda_1)\in {\cal N}$.\newl Let $F_G:\GR\times\GR^2\to
\GC$ be defined by
\begin{equation}\label{FG}
F_G(t,\beta,\lambda_1)=e^{it}\Big[-iv+\beta \sum_{m\in {\mathbb
Z}}\,\frac{g_m(\beta)e^{im{\jmath^*} t}}{i(m{\jmath^*}+1)}+
U(\beta,\lambda_1)(t)\Big].
\end{equation}
Then
${\cal Y}(U(\beta,\lambda_1))(t)=\lambda_1
H\left(F_G(t,\beta,\lambda_1)e^{-it},\mbox{c.c.},\lambda_1\right),$
and, upon setting $z=p-F_G$, (\ref{basic12}) rewrites as
\begin{align*}
\begin{split}
\dot{z}=\lambda_1 e^{it}\big[
H((z+F_G(t,\beta,\lambda_1))e^{-it},\mbox{c.c.},\lambda_1)-H(F_G(t,\beta,\lambda_1)e^{-it},\mbox{\mbox{c.c.}},\lambda_1))\big],
\end{split}
\end{align*}
which reduces to
\begin{equation}
\dot{z}=\lambda_1
e^{it}\widehat{H}(ze^{-it},\ov{z}e^{it},t,\beta,\lambda_1),
\label{newzdotforced}
\end{equation}
where
\begin{equation}\label{thesecondH}\widehat{H}(w,\ov{w},t,\beta,\lambda_1)=H(w+F_G(t,\beta,\lambda_1)e^{-it},\mbox{c.c.},\lambda_1)-H(F_G(t,\beta,\lambda_1)e^{-it},\mbox{c.c.},\lambda_1)\end{equation}
is $2\pi/{\jmath^*}-$periodic in $t$. Then, (\ref{newzdotforced})
becomes
\begin{align}\label{ahahahpolavg}
\begin{split}
\dot{\rho}&= \lambda_1 R_*(\rho)+ \lambda_1 R(t,\psi,\rho,\beta,\lambda_1)\\
\dot{\psi}&= 1+\lambda_1\Psi_*(\rho)+\lambda_1
\Psi(t,\psi,\rho,\beta,\lambda_1),
\end{split}
\end{align}
where $R_*(\rho)= \rho\real \left[L_*(\rho^2)\right]$,
$\Psi_*(\rho)=-\imag\left[L_*(\rho^2)\right]$ for some continuous
function $L_*:\GR\to \GC$, $R,\Psi\in
\mathfrak{P}^{2\pi}_t\cap\mathfrak{ P}^{2\pi}_{\psi}$ and $\Psi$ is
not continuous at $\rho=0$.
\begin{thm}\label{ahahahthmtorus} Assume that $R$ and $\Psi$, as defined in $(\ref{ahahahpolavg})$, are $C^1$ on intervals away from $\rho=0$ and that the averaged
equation
\begin{equation}\label{ahahahaveraged}
\dot{\rho}=\epsilon R_*(\rho)
\end{equation}
has an equilibrium $\rho_0>0$ with $D_{\rho}R_*(\rho_0)=\gamma_0\neq
0$. If the parameters are in a (small enough) deleted neighbourhood
$\mathcal{V}^{\jmath^*}$ of the origin, $(\ref{basic12})$ has an
epicyclic manifold ${\cal G}^{{\jmath^*}}_{\beta,\lambda_1}$, with
$[{\cal G}^{{\jmath^*}}_{\beta,\lambda_1}]_{\DD}=0$. Furthermore,
the origin is a center of drifting when $\lambda_1\gamma_0<0$.
\end{thm}
\begin{rem}
\begin{enumerate}
\item The small term $\hat{\epsilon}g_{-1}(0)$ in (\ref{thefirstH}) and the
absence of a corresponding term in (\ref{thesecondH}) are
responsible for the different form of the regions $\mathcal{V}$ (in
theorem \ref{ahahaveraged}) and $\mathcal{V}^{\jmath^*}$ (in theorem
\ref{ahahahaveraged}), as well as for the location of the center of
drifting. \item The remarks made after theorem~\ref{thmnld1} still
hold, when suitably modified. \item There are a lot of similarities
between our analysis and the results obtained during the analysis of
spiral anchoring in \cite{BLM,Bo1}, such as the presence of
wedge-shaped regions or the deleted neighbourhoods, depending on the
nature of $\jmath^*$. In particular, one might hope that the
epicyclic manifolds would possess $\GZ_{{\jmath^*}}-$spatio-temporal
symmetry; however, this is not the case as the averaged system
defined by (\ref{ahahahpolavg}) is generally only
$2\pi/\jmath^*-$periodic in $t$ when $R\equiv 0$ and $\Psi\equiv 0$.
That being said, the epicycles themselves possess this symmetry in
an appropriate co-rotating frame of reference.
\end{enumerate}
\end{rem}


\section{Epicyclic Drifting For General ESB Terms}
Lifting all restrictions on $\beta$ and $\lambda$ in
(\ref{basic12}), and combining the methods of the preceding section,
we obtain the following general epicyclic drifting theorems for
(\ref{basic12}).
\begin{thm} \label{noclueagainconj} Let $n>1$ and  $k\in
\{1,\ldots,n\}$. Given a hyperbolic equilibrium $\rho_*$ of an
appropriate averaged equation $\dot{\rho}=\beta\tilde{R}_0^k(\rho)$
with eigenvalue $\gamma_k^*$ (derived as in section 3), there is a
region in parameter space near
$(\beta=\lambda_0,\lambda_1,\ldots,\lambda_n)=0$ of the form
$$
{\mathcal V}^1_k=\{(\lambda_0,\lambda_1,\ldots,\lambda_n)\in
\GR^{n+1}\,:\,|\lambda_j|<V_{j,k}|\lambda_k|,\,\,\,V_{j,k}>0,\,\,\mbox{\rm
for $j\neq k$ and $\lambda_k$ near}\,\,0\,\}
$$
when ${\jmath^*}=1$, or of the form \small $$ {\mathcal
V}^{\jmath^*}_k=\{(\beta,\lambda_1,\ldots,\lambda_n)\in
\GR^{n+1}\,:\,|\beta|<\beta_0,\,\,\,|\lambda_j|<V_{j,k}|\lambda_k|,\,\,\,V_{j,k}>0,\,\,\mbox{\rm
for $j\neq k$ and $\lambda_k$ near}\,\,0\,\}
$$ \normalsize
for some constant $\beta_0>0$, when $\jmath^*>1$, such that for all
$0\neq (\beta,\lambda)\in {\mathcal V}_{k}^{\jmath^*}$, with the
additional condition that $\beta\neq 0$ when $\jmath^*=1$,
$(\ref{basic12})$ has an epicyclic manifold ${\cal
E}^{\jmath^*;k}_{\beta,\lambda}$, with $[{\cal
E}^{\jmath^*;k}_{\beta,\lambda}]_{\DD}$ near, but generically not
at, $\xi_k$. Furthermore, $[{\cal
E}^{\jmath^*;k}_{\beta,\lambda}]_{\DD}$ is a center of drifting when
$\lambda_k\gamma_k^*<0$.
\end{thm}
Since the hypotheses of this theorems are not generic, it is not
clear that such integral epicyclic manifolds are common, and their
existence must sometimes be inferred in the physical space,
especially if they are repelling, such as appears to be the case in
\cite{MPMPV}.
\section{An Epicyclic Manifold in Physical Space} In this section, we provide what we believe to be the first observed example of epicyclic drifting in a modified bidomain
system. Our system is a TSB perturbation of the bidomain equations
and parameter values found in \cite{R1}:
\begin{align}\label{thebid}\begin{split} u_t &= \frac{1}{\varsigma} \left(u-\frac{u^3}{3}-v\right)+ \Delta
u+\frac{\alpha \varepsilon
}{1+\alpha(1-\varepsilon)}\psi_{xx} \\
u_{yy}&= \left[\left(2+\alpha+\frac{1}{\alpha}\right)\psi_{xx}+\left(2+\alpha(1-\varepsilon)+\frac{1}{\alpha(1-\varepsilon)}\right)\psi_{yy}\right]\left(1+\frac{1}{\alpha(1-\varepsilon)}\right)^{-1}\frac{1}{\varepsilon} \\
v_t &= \varsigma(u+\delta-\gamma v)+\phi(x-35,y-35)
\end{split}\end{align} where $\varsigma = 0.3$, $\alpha = 1.0$,
$\varepsilon = 0.75$, $\delta = 0.8$, $\gamma = 0.5$ and
$$\phi(z_1,z_2)= -0.03\exp\left(-0.085(z_1^2+z_2^2)\right).$$
The TSB term $\phi(x-35,y-35)$ is uniformly bounded and goes to $0$
as $\|(x,y)\|\to\infty$. Futhermore, it preserves rotations around
the point $(35,35)$. \par As our emphasis lies with qualitative
observations rather than with precise numerical analysis, the
numerical perspective is somewhat naive. The computations are
carried out on a two-dimensional square domain $[10,60]^2$ with 120
grid points to a side and Neumann boundary condition, using a
5-point Laplacian, continuous linear finite elements on square
meshes and the fully implicit second order Gear finite difference
integrator. The tip path of the $u$ component of two solutions are
shown in  \begin{figure}[t]
\begin{center}
\includegraphics[width=222pt]{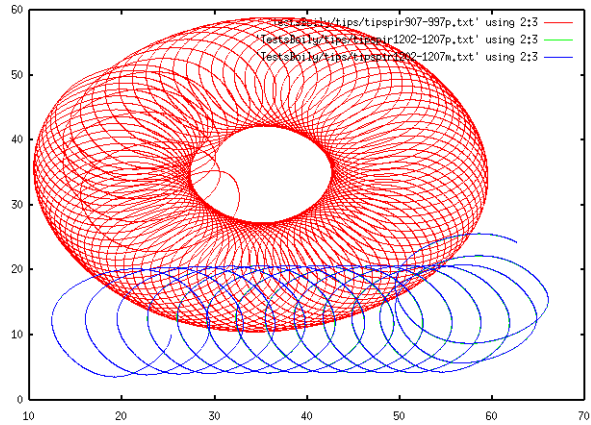}
\caption{Tip path of the $u$ component of two solutions of
(\ref{thebid}). The solution in red is attracted by an epicyclic
manifold. The solution in blue shows the effect of the boundary: the
two solutions are clearly not of the same nature. Compare the
epicyclic manifold with the image in figure
\ref{epmannn}.}\end{center}\hrule
\end{figure}
figure~\thefigure.

\section{Summary and Concluding Remarks}
{Recently, equivariant dynamical systems theory has been used to
provide an approach to the study of spiral wave dynamics and
bifurcations, in particular, it has provided mechanisms for such
behaviour as spiral tip meandering and resonant growth
\cite{B2,SSW1,SSW2,SSW3,SSW4,Wulff}, spiral anchoring/repelling and
boundary drifting \cite{LW}, which have been explained as
consequences of forced Euclidean symmetry breaking.}

{In this paper, we have used this model-independent approach to
analyze epicyclic drifting of the spiral tip in media with several
localized inhomogeneities, with or without anisotropy. The result of
a simple numerical experiment is in agreement with our theoretical
conclusions. The RSB terms are characterized by an integer
$\jmath^*>1$. It is important to note that, as of now, only the
integers $\jmath^*=2$ (anisotropic cardiac tissue, say) and
$\jmath^*=1$ (an excitable medium that in which there is a directed
current, such as a reaction-diffusion-advection system) have easy
interpretations in the context of RSB.}

{It should be recalled that our analysis rests on certain
simplifying modeling assumptions which may not be valid in some
realistic physical systems: namely concerning the discrete nature of
the inhomogeneities (i.e. finite number of inhomogeneity sites) and
the hypothesis of local rotational symmetry of the individual
inhomogeneities. In a realistic model of excitable media such as the
bidomain model describing electrical conduction of cardiac tissue,
with or without advection, actual inhomogeneities may lack this
local rotational symmetry, or may even be distributed smoothly and
non-symmetrically over the medium (\textit{i.e.} an inhomogeneity
field).} \par Should the discrete localized inhomogeneities not
possess circular symmetry, we believe that our results would still
describe the essential qualitative features of epicyclic drifting,
even though our analysis does not technically apply, as epicyclic
drifting is linked to hyperbolic fixed points. However, for a
smoothly distributed non-symmetric inhomogeneity field, our
techniques are unlikely to yield meaningful results.\newl Finally,
we would like to point out that our existence results do not provide
a description of the flow on the epicyclic manifolds: we gave two
examples of center bundle systems with essentially different flows
in \cite{Byeah}. In the first example, the flow is ``ergodic'' (see
the first figure of this article): if the spiral tip lies in the
epicyclic manifold, the tip path eventually fills the entire
manifold. In the second example, the flow on the epicyclic manifold
is dictated by the stability of rotating waves located \textit{on}
the manifold (see
\begin{figure}[t]
\begin{center}
\includegraphics[width=240pt]{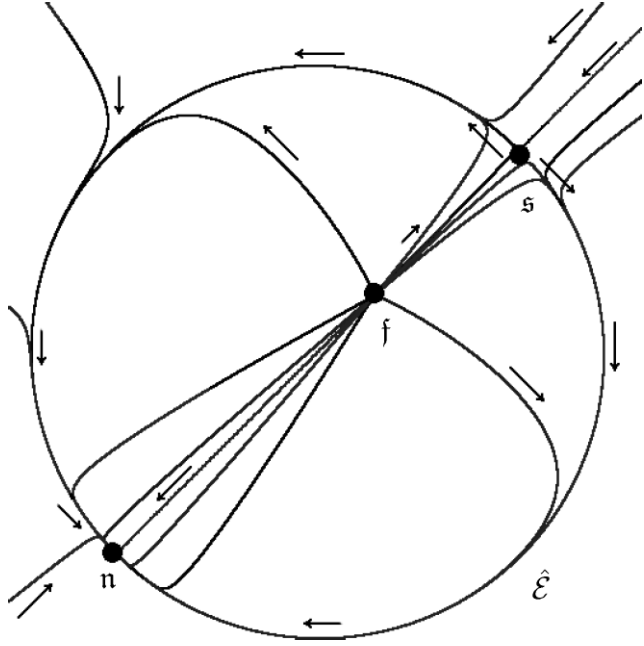}
\caption{Projection (in the $z-$plane) of a stable epicycle manifold
$\hat{\cal E}$ containing three rotating waves $\mathfrak{s}$,
$\mathfrak{n}$ and $\mathfrak{f}$. The arrows indicate the flow on
the manifold. \label{epmannnnn}}\end{center}\hrule
\end{figure}
figure~\thefigure).

\section{Acknowledgements}
{The author is grateful to the referees for suggestions that have
made this paper much more readable. The author would also like to
thank Prof. Victor G. LeBlanc for the helpful advice and feedback he
has provided over the last 10 years, and Marc Ethier and Eric Matsui
for lending a hand with some of the numerics.}

\appendix
\section{Technical Results} \textbf{Proposition A.1}
\textit{The function $M^1(w,\ov{w},0)$ defined in $(\ref{avgM})$ is
$\GS^1-$equivariant.}\newl\small \proof Recall that $\hat{H}_1$ is
$\GS^1-$equivariant by construction. Then
\begin{align*}
M^1(we^{-i\theta},\ov{w}e^{i\theta},0)&=\frac{1}{2\pi}\int_{0}^{2\pi}e^{it}K(we^{-i\theta}e^{-it},\ov{w}e^{i\theta}e^{it},t,0,0)\, dt \\
&=\frac{1}{2\pi}\int_{0}^{2\pi}e^{it}\hat{H}_1(we^{-it}e^{-i\theta},\ov{w}e^{it}e^{i\theta},0,0)\, dt \\
&=\frac{1}{2\pi}\int_{0}^{2\pi}e^{it}e^{i\theta}\hat{H}_1(we^{-it},\ov{w}e^{it},0,0)\, dt \\
&=e^{i\theta}M^1(w,\ov{w},0),
\end{align*} that is, $M^1$ is $\GS^1-$equivariant. \qed\newl\normalsize
\textbf{Proposition A.2} \textit{Let all terms, variables and
functions be as in Theorem~\ref{thmtorus}. \label{appH1H2a} In
particular, the functions $R$ and $\Psi$ are $C^1$ on intervals away
from $\rho=0$. Denote
\begin{align*}
\Sigma^r_0&=[0,2\pi]\times[0,2\pi]\times\{0\}\times S_0 \\
\Sigma^r_{\sigma}&=[0,2\pi]\times[0,2\pi]\times [-\sigma,\sigma]\times S_0.
\end{align*} Note that these spaces, as well as the spaces $\Sigma_0$ and
$\Sigma_{\sigma}$ from Theorem~\ref{Halesavg}, are convex. The
functions $\Theta$ and $\Lambda$ satisfy the following conditions:}
\begin{enumerate}
\item \textit{$\Theta$ and $\Lambda$ are bounded by a function $\Xi(\epsilon,\nu)=O(\epsilon,\nu_2,\ldots,\nu_n)$ over $\Sigma_0,$ and}
\item \textit{for all $0\leq \sigma\leq \sigma_0$, $\Theta$ and $\Lambda$ are Lipschitz in Hale's sense (with Lipschitz constants $\theta(\epsilon,\nu,\sigma)=O(\epsilon,\nu_2,\ldots,\nu_n,\sigma)$ and $\nu(\epsilon,\nu,\sigma)=O(\epsilon,\nu_2,\ldots,\nu_n,\sigma)$, respectively) over $\Sigma_{\sigma}$.} (p.~\pageref{pappH1H2})
\end{enumerate}
\small \proof Since $\Theta,\Lambda\in \mathfrak{P}^{2\pi}_t\cap
\mathfrak{P}^{2\pi}_{\psi}$, we need only show that the first statement holds
over $\Sigma^r_0$ and the second over $\Sigma^r_{\sigma}$, for all $0\leq
\sigma\leq \sigma_0$. For $j=1, \ldots, n$, there are appropriate functions
$R_j,\Psi_j\in \mathfrak{P}^{2\pi}_t\cap \mathfrak{P}^{2\pi}_{\psi}$, $C^1$ on
intervals away from $\rho=0$, such that
\begin{align}
\label{thedefin}
\begin{split}
\Psi(t,\psi,\rho,\epsilon,\nu)&=\epsilon \Psi_1(t,\psi,\rho,\epsilon,\nu)+\sum_{j=2}^n\nu_j\Psi_j(t,\psi,\rho,\nu) \\
R(t,\psi,\rho,\epsilon,\nu)&=\epsilon
R_1(t,\psi,\rho,\epsilon,\nu)+\sum_{j=2}^n\nu_jR_j(t,\psi,\rho,\nu),
\end{split}
\end{align}
according to (\ref{funcW}).
\begin{enumerate}
\item Over $\Sigma^r_0$, we have $x=0$ and so \begin{align}\label{upthere1}\begin{split}
\Theta(t,\psi,0,\epsilon,\nu)&=\epsilon\Psi(t,\psi,\rho(\epsilon,\nu),\epsilon,\nu)\\&
=\epsilon^2\Psi_1(t,\psi,\rho(\epsilon,\nu),\epsilon,\nu)+\epsilon
\sum_{j=2}^n\nu_j\Psi_j(t,\psi,\rho(\epsilon,\nu),\epsilon,\nu)
\end{split}\end{align}
and
\begin{align}\label{upthere2}\begin{split}
\Lambda(t,\psi,0,\epsilon,\nu)&=R(t,\psi,\rho(\epsilon,\nu),\epsilon,\nu)\\&
=\epsilon R_1(t,\psi,\rho(\epsilon,\nu),\epsilon,\nu)+
\sum_{j=2}^n\nu_j
R_j(t,\psi,\rho(\epsilon,\nu),\epsilon,\nu)\end{split}\end{align}
according to (\ref{thedefin}). \par For $j=1,\ldots, n$, the
continuous functions $|R_j|$ and $|\epsilon\Psi_j|$ reach their
maximum $C_j$ and $E_j$, respectively, on the compact set
$[0,2\pi]\times \{0\}$ for $j=1,\ldots, n$. Then
$$\left|\epsilon\Psi_j(t,\psi,\rho(\epsilon,\nu),\epsilon,\nu)\right|\leq E_j\quad\mbox{and}\quad \left|R_j(t,\psi,\rho(\epsilon,\nu),\epsilon,\nu)\right|\leq C_j$$ over $\Sigma_0^r$ for $j=1,\ldots,n$.
According to (\ref{upthere1}) and (\ref{upthere2}),
\begin{align*}
|\Theta(t,\psi,0,\epsilon,\nu)|&\leq |\epsilon| E_1+\sum_{j=2}^n |\nu_j|E_j=Q_1(\epsilon,\nu) \\
|\Lambda(t,\psi,0,\epsilon,\nu)|&\leq |\epsilon| C_1
+\sum_{j=2}^n|\nu_j| C_j=Q_2(\epsilon,\nu)
\end{align*} over $\Sigma^r_0$. Set $$\Xi(\epsilon,\nu)=\max\{Q_1(\epsilon,\nu),Q_2(\epsilon,\nu)\}.$$ Then $\Lambda$ and $\Theta$ are bounded by $\Xi(\epsilon,\nu)=O(\epsilon,\nu_2,\ldots,\nu_n)$ over~$\Sigma^r_0$.
\item \label{appH1H2} Let $(\psi_1,x_1),(\psi_2,x_2)\in \GR\times [-\sigma,\sigma]$ for $0\leq\sigma\leq\sigma_0$. By one of the mean value theorems, there exist points $(\psi^*,x^*),(\psi_*,x_*)\in [0,2\pi]\times [-\sigma,\sigma]$ on the line joining $(\psi_1,x_1)$ and $(\psi_2,x_2)$ such that \begin{align*}
|\Theta(t,\psi_1,x_1,\epsilon,\nu)-\Theta(t,\psi_2,x_2,\epsilon,\nu) |&=|\widehat{\Theta}(t,\psi^*,x^*,\epsilon,\nu)| \big[|\psi_1-\psi_2|+|x_1-x_2|\big] \\
|\Lambda(t,\psi_1,x_1,\epsilon,\nu)-\Lambda(t,\psi_2,x_2,\epsilon,\nu)
|&=|\widehat{\Lambda}(t,\psi_*,x_*,\epsilon,\nu)|\big[|\psi_1-\psi_2|+|x_1-x_2|\big]
,\end{align*}where
\begin{align*} \widehat{\Theta}(t,\psi,x,\epsilon,\nu)&=D_x\Theta(t,\psi,x,\epsilon,\nu)+D_{\psi}\Theta(t,\psi,x,\epsilon,\nu) \\ &=x K_0^{\Psi}(x,\epsilon,\nu)+\epsilon K_1^{\Psi}(t,\psi,x,\epsilon,\nu)+\sum_{j=2}^n\nu_j K_j^{\Psi}(t,\psi,x,\epsilon,\nu) \\
\widehat{\Lambda}(t,\psi,x,\epsilon,\nu)&=D_x\Lambda(t,\psi,x,\epsilon,\nu)+D_{\psi}\Lambda(t,\psi,x,\epsilon,\nu)
\\ &=x K_0^{R}(x,\epsilon,\nu)+\epsilon
K_1^{R}(t,\psi,x,\epsilon,\nu)+\sum_{j=2}^n\nu_j
K_j^{R}(t,\psi,x,\epsilon,\nu), \end{align*} where
\begin{align*}K_0^{\Psi}(x,\epsilon,\nu)&=\epsilon D_x B_2(x,\epsilon,\nu) \\ K_0^{R}(x,\epsilon,\nu)&=\epsilon \left(D_x B_1(x,\epsilon,\nu)x + 2B_1(x,\epsilon,\nu)\right)\\
 K_1^{\Psi}(t,\psi,x,\epsilon,\nu)&=B_2(x,\epsilon,\nu)+\epsilon\left(D_x \Psi_1(t,\psi,\rho(\epsilon,\nu)+x,\epsilon,\nu)+D_{\psi} \Psi_1(t,\psi,\rho(\epsilon,\nu)+x,\epsilon,\nu)\right) \\ K_1^{R}(t,\psi,x,\epsilon,\nu)&=\epsilon\left(D_x R_1(t,\psi,\rho(\epsilon,\nu)+x,\epsilon,\nu)+D_{\psi} R_1(t,\psi,\rho(\epsilon,\nu)+x,\epsilon,\nu)\right) \\ \intertext{and}
 K_j^{\Psi}(t,\psi,x,\epsilon,\nu)&=\epsilon\left(D_x \Psi_j(t,\psi,\rho(\epsilon,\nu)+x,\epsilon,\nu)+D_{\psi} \Psi_j(t,\psi,\rho(\epsilon,\nu)+x,\epsilon,\nu)\right) \\ K_j^{R}(t,\psi,x,\epsilon,\nu)&=\epsilon\left(D_x R_j(t,\psi,\rho(\epsilon,\nu)+x,\epsilon,\nu)+D_{\psi} R_j(t,\psi,\rho(\epsilon,\nu)+x,\epsilon,\nu)\right),
\end{align*} where $\Psi_j$ and $R_j$ are as in (\ref{thedefin}).
Since $\Theta$ and $\Lambda$ are continuously differentiable,
$\widehat{\Theta}$ and $\widehat{\Lambda}$ are continuous on
$\Sigma^r_{\sigma}$, as are $K_j^{\Psi}$ and $K_j^R$ for $j=0,\ldots, n$. \par
In particular, the functions $|K_j^{\Psi}|$ and $|K_j^{\Psi}|$ each reach their
respective maximum $k_j^{\Psi}$ and $k_j^{R}$ on $\Sigma^r_{\sigma}$ for
$j=0,\ldots, n$. Then, note that $|x^*|,|x_*|\leq \sigma$,
\begin{align*}|\widehat{\Theta}(t,\psi^*,x^*,\epsilon,\nu)|&\leq
|x^*| |K_0^{\Psi}(x^*,\epsilon,\nu)|+|\epsilon|
|K_1^{\Psi}(t,\psi^*,x^*,\epsilon,\nu)|+\sum_{j=2}^n|\nu_j|
|K_j^{\Psi}(t,\psi^*,x^*,\epsilon,\nu)| \\ &\leq |x^*|
k_0^{\Psi}+|\epsilon| k_1^{\Psi}+\sum_{j=2}^n |\nu_j| k_j^{\Psi}
\\ &\leq \sigma k_0^{\Psi}+|\epsilon| k_1^{\Psi}+\sum_{j=2}^n
|\nu_j| k_j^{\Psi}=\theta(\epsilon,\nu,\sigma)
\end{align*}
and
\begin{align*}|\widehat{\Lambda}(t,\psi_*,x_*,\epsilon,\nu)|&\leq
|x_*| |K_0^{R}(x_*,\epsilon,\nu)|+|\epsilon|
|K_1^{R}(t,\psi_*,x_*,\epsilon,\nu)|+\sum_{j=2}^n|\nu_j|
|K_j^{R}(t,\psi_*,x_*,\epsilon,\nu)| \\ &\leq |x_*|
k_0^{R}+|\epsilon| k_1^{R}+\sum_{j=2}^n |\nu_j| k_j^{R} \\ &\leq
\sigma k_0^{R}+|\epsilon| k_1^{R}+\sum_{j=2}^n |\nu_j|
k_j^{R}=\nu(\epsilon,\nu,\sigma).
\end{align*} In particular \begin{align*}|\Theta(t,\psi_1,x_1,\epsilon,\nu)-\Theta(t,\psi_2,x_2,\epsilon,\nu) |&\leq \theta(\epsilon,\nu_2,\ldots,\nu_n,\sigma)\big[|\psi_1-\psi_2|+|x_1-x_2|\big] \\ |\Lambda(t,\psi_1,x_1,\epsilon,\nu)-\Lambda(t,\psi_2,x_2,\epsilon,\nu) |&\leq \nu(\epsilon,\nu_2,\ldots,\nu_n,\sigma)\big[|\psi_1-\psi_2|+|x_1-x_2|\big], \end{align*} where $\theta(\epsilon,\nu,\sigma),\nu(\epsilon,\nu,\sigma)=O(\epsilon,\nu_2,\ldots,\nu_n,\sigma)$. Hence $\Theta$ and $\Lambda$ are Lipschitz in Hale's sense.
\end{enumerate} This completes the proof. \qed\normalsize
\bibliographystyle{amsplain}
\bibliography{spiralsp}

\providecommand{\BIBYu}{Yu}
\providecommand{\bysame}{\leavevmode\hbox to3em{\hrulefill}\thinspace}
\providecommand{\MR}{\relax\ifhmode\unskip\space\fi MR }
\providecommand{\MRhref}[2]{%
  \href{http://www.ams.org/mathscinet-getitem?mr=#1}{#2}
}
\providecommand{\href}[2]{#2}
\begin{thebibliography}{10}

\bibitem{B1}
D.~Barkley, \emph{Linear stability analysis of rotating spiral waves in
  excitable media}, Phys. Rev. Lett. \textbf{68} (1992), 2090--3.

\bibitem{B2}
\bysame, \emph{Euclidean symmetry and the dynamics of rotating spiral waves},
  Phys. Rev. Lett. \textbf{76} (1994), 164--7.

\bibitem{BK}
D.~Barkley and I.~G. Kevrekedis, \emph{A dynamical system approach to spiral
  wave dynamics}, Chaos \textbf{4} (1994), 453--60.

\bibitem{BKT}
D.~Barkley, M.~Kness, and L.~S. Tuckerman, \emph{Spiral-wave dynamics in a
  simple model of excitable media: The transition from simple to compound
  rotation}, Phys. Rev. Lett. \textbf{42} (1990), 2489--92.

\bibitem{Byeah}
P.~Boily, \emph{Spiral wave dynamics under full euclidean symmetry-breaking: A
  dynamical system approach}, Ph.D. thesis, University of Ottawa, 2006.

\bibitem{Bo1}
\bysame, \emph{Spiral anchoring in anisotropic media with inhomogeneities},
  ([arXived]).

\bibitem{BLM}
P.~Boily, V.~G LeBlanc, and E.~Matsui, \emph{Spiral anchoring in media with
  multiple inhomogeneities: a dynamical system approach}, J. Nonlin. Sc. ([to
  be published]).

\bibitem{BEL}
Y.~Bourgault, M.~Ethier, and V.G. LeBlanc, \emph{Simulation of
  electrophysiological waves with an unstructured finite element method},
  ESAIM:M2AN \textbf{37} (2003), 649--62.

\bibitem{Detal}
J.~M. Davidenko, A.~V. Persov, R.~Salomonsz, W.~Baxter, and J.~Jalife,
  \emph{Stationary and drifting spiral waves of excitation in isolated cardiac
  muscle}, Nature \textbf{355} (1992), 349--51.

\bibitem{DMcK}
H.~Dym and H.~P. McKean, \emph{Fourier series and integrals}, Academic Press,
  New York, 1972.

\bibitem{FSSW}
B.~Fiedler, B.~Sandstede, A.~Scheel, and C.~Wulff, \emph{Bifurcation from
  relative equilibria of noncompact group actions: Skew products, meanders and
  drifts}, Doc. Math. \textbf{1} (1996), 479--555.

\bibitem{GSS}
M.~Golubitsky and D.~G. Stewart, I. et~Schaeffer, \emph{Singularities and
  groups in bifurcation theory, volume ii}, Springer-Verlag, Berlin, 1988.

\bibitem{GZM}
S.~Grill, V.~S. Zykov, and S.~C. M{\"{u}}ller, \emph{Spiral wave dynamics under
  pulsatory modulation of excitability}, J. Phys. Chem. \textbf{100} (1996),
  19082--8.

\bibitem{Grubb}
N.R. Grubb and S.~Furniss, \emph{Science, medicine and the future:
  Radiofrequency ablation for atrial fibrillation}, British Medical Journal
  \textbf{322} (2001), 777--780.

\bibitem{H1}
J.~K. Hale, \emph{Integral manifolds of perturbed differential equations}, Ann.
  Math. \textbf{73} (1961), 496--531.

\bibitem{J}
J.~Jalife, \emph{Ventricular fibrillation: Mechanisms of initiation and
  maintenance}, Annu. Rev. Physiol. \textbf{62} (2000), 25--50.

\bibitem{KS}
J.~Keener and J.~Sneyd, \emph{Mathematical physiology}, IAM, Springer, New
  York, 1998.

\bibitem{LeB}
V.~G. LeBlanc, \emph{Rotational symmetry-breaking for spiral waves},
  Nonlinearity \textbf{15} (2002), 1179--203.

\bibitem{LW}
V.~G. LeBlanc and C.~Wulff, \emph{Translational symmetry-breaking for spiral
  waves}, J. Nonlin. Sc. \textbf{10} (2000), 569--601.

\bibitem{LR}
V.G. LeBlanc and B.J. Roth, \emph{Meandering of spiral waves in anisotropic
  tissue}, Dynam. Contin. Discrete Impuls. Systems \textbf{B10} (2003), 29--42.

\bibitem{LOPS}
G.~Li, Q.~Ouyang, V.~Petrov, and H.~L. Swinney, \emph{Transition from simple
  rotating chemical spirals to meandering and traveling spirals}, Phys. Rev.
  Lett. \textbf{77} (1996), 2105--9.

\bibitem{MAK}
D.~Mackenzie, \emph{Making sense of a heart gone wild}, Science \textbf{303}
  (2004), 786--7.

\bibitem{MDC}
J.E. Marine, J.~Dong, and H.~Calkins, \emph{Catheter ablation therapy for
  atrial fibrillation}, Progress in Cardiovascular Diseases \textbf{48} (2005),
  no.~3, 178--192.

\bibitem{MDZ}
A.S Mikhailov, V.A. Davydov, and V.S. Zykov, \emph{Complex dynamics of spiral
  waves and motion of curves}, Physica D \textbf{70} (1994), 1--39.

\bibitem{MPMPV}
A.~P. Mu{\~{n}}uzuri, V.~P{\'{e}}rez-Mu{\~{n}}uzuri, and V.~P{\'e}rez-Villar,
  \emph{Attraction and repulsion of spiral waves by localized inhomogeneities
  in excitable media}, Phys. Rev. E \textbf{58} (1998), R2689--92.

\bibitem{R2}
B.~J. Roth, \emph{Approximate analytical solutions to the bidomain equations
  with unequal anisotropy ratios}, Phys. Rev. E \textbf{55} (1997), 1819--26.

\bibitem{R1}
\bysame, \emph{Frequency locking of meandering spiral waves in cardiac tissue},
  Phys. Rev. E \textbf{57} (1998), R3735--8.

\bibitem{abscon2}
B.~Sandstede and A.~Scheel, \emph{Absolute and convective instabilities of
  waves on unbounded and large bounded domains}, Physica D \textbf{145} (2000),
  233--277.

\bibitem{abscon1}
\bysame, \emph{Absolute versus convective instability of spiral waves}, Phys.
  Rev. E (3) \textbf{62} (2000), 7708--7714.

\bibitem{SSW1}
B.~Sandstede, A.~Scheel, and C.~Wulff, \emph{Center manifold reduction for
  spiral waves}, C. R. Acad. Sci. \textbf{324} (1997), 153--8.

\bibitem{SSW2}
\bysame, \emph{Dynamics of spiral waves on unbounded domains using
  center-manifold reductions}, J. Diff. Eq. \textbf{141} (1997), 122--49.

\bibitem{SSW3}
\bysame, \emph{Bifurcations and dynamics of spiral waves}, J. Nonlin. Sc.
  \textbf{9} (1999), 439--78.

\bibitem{SSW4}
\bysame, \emph{Dynamical behavior of patterns with euclidean symmetry}, Pattern
  Formation in Continuous and Coupled Systems (Golubitsky, Luss, and Strogatz,
  eds.), Springer-Verlag, Berlin, 1999.

\bibitem{S}
A.~Scheel, \emph{Bifurcation to spiral waves in reaction-diffusion systems},
  SIAM J. Math. Anal. \textbf{29} (1998), 1399--418.

\bibitem{WR}
N.~Wiener and A.~Rosenblueth, \emph{The mathematical formulation of the problem
  of conduction of impulses in a network of connected excitable elements,
  specifically in cardiac muscle}, Arch. Inst. Card. De Mexico \textbf{16}
  (1946), 205--65.

\bibitem{W}
A.~T. Winfree, Cardiac Electrophysiology, From Cell to Bedside (Zipes and
  Jalife, eds.), Saunders, Philadelphia, second ed., 1995, pp.~379--89.

\bibitem{Wetal}
F.~X. Witkowski, L.~J. Leon, P.~A. Penkoske, W.~R. Giles, Mark~L. Spanol, W.~L.
  Ditto, and A.~T. Winfree, \emph{Spatiotemporal evolution of ventricular
  fibrillation}, Nature \textbf{392} (1998), 78--82.

\bibitem{Wulff}
C.~Wulff, \emph{Theory of meandering and drifting spiral waves in
  reaction-diffusion systems}, Ph.D. thesis, Freie Universit{\"{a}}t Berlin,
  1996.

\bibitem{YP}
Y.~A. Yermakova and A.~M. Pertsov, \emph{Interaction of rotating spiral waves
  with a boundary}, Biophys. \textbf{31} (1987), 932--40.

\bibitem{ZM}
V.~S. Zykov and S.~C. M{\"{u}}ller, \emph{Spiral waves on circular and
  spherical domains of excitable medium}, Physica D \textbf{97} (1996),
  322--32.

\end{thebibliography}
\end{document}